\documentclass[11pt,reqno]{article}

\usepackage{geometry}

\usepackage[pdftitle={Holomorphic Poisson Structures and Groupoids},
pdfauthor={Camille Laurent, Mathieu Stienon, Ping Xu}]{hyperref}

\usepackage[T1]{fontenc}
\usepackage[latin1]{inputenc}
\usepackage{lmodern}

\usepackage{setspace}
\usepackage{indentfirst}

\usepackage[intlimits]{amsmath}
\usepackage{amssymb}
\usepackage{stmaryrd}

\usepackage{graphicx}
\usepackage[all]{xy}
\def\dar[#1]{\ar@<2pt>[#1]\ar@<-2pt>[#1]}
\usepackage{amsthm}

\theoremstyle{plain}% default
\newtheorem{prop}{Proposition}[section]
\newtheorem{lem}[prop]{Lemma}

\newtheorem{thm}[prop]{Theorem}

\newtheorem*{prop*}{Proposition}
\newtheorem*{lem*}{Lemma}
\newtheorem*{sublem*}{Sublemma}
\newtheorem*{cor*}{Corollary}
\newtheorem*{thm*}{Theorem}
\newtheorem*{hypo*}{Hypothesis}
\newtheorem*{question*}{Question}
\newtheorem*{conjecture*}{Conjecture}
\newtheorem*{scholum*}{Scholum}
\newtheorem{defn}[prop]{Definition}
\newtheorem*{defn*}{Definition}

\theoremstyle{definition}

\newtheorem*{con*}{Construction}
\newtheorem*{note*}{Note}

\theoremstyle{remark}

\newtheorem*{warning*}{Warning}
\newtheorem*{shortnote*}{Note}
\newtheorem*{claim*}{Claim}
\newtheorem*{axiom*}{Axiom}

\newtheoremstyle{slanted}% name
  {3pt}%      Space above, empty = `usual value'
  {3pt}%      Space below
  {\slshape}% Body font
  {}%         Indent amount (empty = no indent, \parindent = para indent)
  {\bfseries}% Thm head font
  {.}%        Punctuation after thm head
  {.5em}%     Space after thm head: " " = normal interword space;
  {}% Thm head spec

\theoremstyle{slanted}
\newtheorem{example}[prop]{Example}

\newtheorem*{example*}{Example}
\newtheorem*{examples*}{Examples}

\newtheorem*{ex*}{Example}
\newtheorem*{exs*}{Examples}
\newtheorem{remark}[prop]{Remark}

\newtheorem*{remark*}{Remark}
\newtheorem*{remarks*}{Remarks}
\newtheorem{rmk}[prop]{Remark}

\newtheorem*{rmk*}{Remark}
\newtheorem*{rmks*}{Remarks}

\DeclareMathOperator{\idn}{id}

\DeclareMathOperator{\Lie}{Lie}

\newcommand{\beq}[1]{\begin{equation}\label{#1}}
\newcommand{\eeq}{\end{equation}}

\newcommand{\CC}{\mathbb{C}}
\newcommand{\RR}{\mathbb{R}}

\newcommand{\lie}[2]{[#1,#2]} % Lie bracket
\newcommand{\schouten}[2]{[#1,#2]} % Schouten bracket
\newcommand{\delbar}{\bar{\partial}}

\newcommand{\half}{\frac{1}{2}}
\newcommand{\rond}{\circ}

\newcommand{\genrel}[2]{\left\{ #1 | #2 \right\}}

\newcommand{\cinf}[1]{C^{\infty}(#1)}
\newcommand{\sections}[1]{\Gamma(#1)}
\newcommand{\vf}{\mathfrak{X}} % vector fields
\newcommand{\toto}{\rightrightarrows}

\newcommand{\diese}{^{\sharp}}

\newcommand{\inv}{^{-1}}

\newcommand{\ddtz}[1]{\left.\tfrac{d}{dt}#1\right|_0}

\newcommand{\LLie}{\mathbb{L}\text{\upshape{ie}}}
\newcommand{\NIJ}{\phi} % de Rham
\newcommand{\diff}{{\rm d}} % de Rham
\newcommand{\gm}{\Gamma}
\newcommand{\product}{m}
\newcommand{\unit}{\epsilon}
\newcommand{\inverse}{\iota}
\newcommand{\hf}[1]{\mathcal{O}_{#1}} % sheaf of holomorphic functions on
\newcommand{\shs}[1]{\mathcal{#1}} % stands for sheaf of holomorphic sections of holomorphic bundle
\newcommand{\pire}{\pi_R} % real part of bivector \pi
\newcommand{\piim}{\pi_I} % imaginary part of bivector \pi
\newcommand{\pb}[2]{\{#1,#2\}} % Poisson bracket on functions associated to holomorphic \pi
\newcommand{\pbre}[2]{\{#1,#2\}_R} % Poisson bracket on functions associated to \pire
\newcommand{\pbim}[2]{\{#1,#2\}_I} % Poisson bracket on functions associated to \piim
\newcommand{\be }{\begin{eqnarray*}}
\newcommand{\ee }{\end{eqnarray*}}

\newcommand{\lon }{\longrightarrow }

\newcommand{\widetild}[1]{\mathcal{S}#1}

\begin{document}

%\title[]{Holomorphic Poisson Structures and Groupoids}
%\author[]{Camille Laurent-Gengoux}
%\address{D\'epartement de math\'ematiques, Universit\'e de Poitiers, 
%86962 Futuroscope-Chasseneuil, France}
%\email{laurent@math.univ-poitiers.fr}
%\author[]{Mathieu Sti\'enon}
%\thanks{Research supported by the European Union through the FP6 Marie Curie R.T.N. ENIGMA
%(Contract number MRTN-CT-2004-5652).}
%\address{E.T.H.~Z\"urich, Departement Mathematik, 8092 Z\"urich, Switzerland}
%\email{stienon@math.ethz.ch}
%\author[]{Ping Xu}
%\address{Department of Mathematics, Penn State University, University Park, PA 16802, U.S.A.}
%\email{ping@math.psu.edu}
%\date{\texttt{\jobname.tex}}
%%\begin{abstract}
%%\end{abstract}
%%\subjclass{}
%\maketitle

\title{Integration of Holomorphic Lie Algebroids}
\author{
Camille Laurent-Gengoux \\ 
D\'epartement de math\'ematiques \\ Universit\'e de Poitiers \\ 
86962 Futuroscope-Chasseneuil, France \\ 
\href{mailto:laurent@math.univ-poitiers.fr}{\texttt{laurent@math.univ-poitiers.fr}}
\and 
Mathieu Sti\'enon
\thanks{Research supported by the European Union through the FP6 Marie Curie R.T.N. ENIGMA
(Contract number MRTN-CT-2004-5652).} \\
E.T.H.~Z\"urich \\ Departement Mathematik \\ 
8092 Z\"urich, Switzerland \\
\href{mailto:stienon@math.ethz.ch}{\texttt{stienon@math.ethz.ch}}
\and 
Ping Xu
\thanks{Research partially supported by NSF
grants DMS-0306665 and DMS-0605725 \&  NSA grant H98230-06-1-0047} \\
Department of Mathematics \\ Penn State University \\ 
University Park, PA 16802, U.S.A. \\
\href{mailto:ping@math.psu.edu}{\texttt{ping@math.psu.edu}} }
%\date{\texttt{\jobname.tex}}
\date{}
%\subjclass{}
\maketitle

%\centerline{{Dedicated to the memory of Paulette Libermann}}

\begin{abstract}
We prove that a holomorphic Lie algebroid is
integrable if, and only if, its underlying real Lie
algebroid is integrable. Thus the integrability criteria of Crainic-Fernandes
(Theorem 4.1 in \cite{CF})
do also apply in the holomorphic context without any modification.
As a consequence we gave another  proof of the theorem
 that a holomorphic Poisson manifold is integrable if and only if 
its real part or  imaginary part is
 integrable as a  real Poisson manifold \cite{crainic, PqN}.
\end{abstract}

\tableofcontents

\section{Introduction}

Since Lie's third theorem fails for  Lie algebroids,
it has been a central theme of study in the theory of
Lie groupoids whether a given   Lie algebroid
is integrable. By an integrable Lie algebroid, we
mean there exists an $s$-connected and $s$-simply connected
 Lie groupoid of which it is the infinitesimal version.
For real Lie algebroids, the integrability  problem
has been completely solved by Crainic-Fernandes \cite{CF}
based on the work of  Cattaneo-Felder on Poisson
sigma models \cite{CattaneoF}.

Recently, there has been increasing interest
in  holomorphic Lie algebroids and  holomorphic
Lie groupoids.  It is a very natural question
to find the  integrability condition for 
holomorphic Lie algebroids.
In \cite{LSX}, we systematically
 studied holomorphic Lie algebroids and their
relation with real Lie algebroids. In particular, we proved that
associated to any holomorphic Lie algebroid $A$, there is
a  canonical real  Lie algebroid $A_R$ such that
the inclusion $\shs{A}\to \shs{A}_\infty$
is a morphism of sheaves.
Here $\shs{A}$ and $\shs{A}_\infty$ denote, respectively,
   the sheaf of holomorphic 
sections and the sheaf of smooth sections   of $A$.
In other words, a holomorphic Lie algebroid
can be considered as  a holomorphic vector bundle
 $A\to X$
 whose underlying real vector
bundle is endowed with a Lie algebroid structure
such that, for any open subset $U\subset X$,
$\lie{\shs{A}(U)}{\shs{A}(U)}\subset\shs{A}(U)$ 
 and the restriction of the Lie bracket
$\lie{\cdot}{\cdot}$ to $\shs{A}(U)$ is $\CC$-linear.
It is thus natural to ask

\begin{quote}
\textbf{Problem 1.}
Given a holomorphic Lie algebroid $A$ with underlying real Lie
algebroid $A_R$, what is the relation between the integrability
of $A$ and the integrability of $A_R$?
\end{quote}

To tackle this problem, we need an equivalent  description 
of holomorphic Lie algebroids. For this purpose, we will
view a  holomorphic Lie algebroid
as a real Lie algebroid structure on a
 holomorphic vector bundle $A\to X$,
 whose almost complex structure $J_A :TA\to TA$
is an infinitesimal multiplicative $(1,1)$-tensor.
By {\em an infinitesimal multiplicative $(1,1)$-tensor}
on a (real)
Lie algebroid $A$, we mean a $(1,1)$-tensor $\NIJ_A$ on $A$
such that $\NIJ_A:TA\to TA$ is an  endomorphism
of the tangent Lie algebroid $TA\to TX$.

 On the level of Lie groupoids, if $\gm\toto X$ is an $s$-connected and $s$-simply connected Lie groupoid integrating the real Lie algebroid $A$, then $J_A$ integrates to an automorphism $J_\gm: T\gm \to T\gm$ of the tangent groupoid $T\gm\toto TX$ such that $J_\gm^2=-\idn$. 
One can show that $J_\gm$ is fiberwise linear with respect to the vector
bundle structure $T\gm\to\gm$. Hence it is a $(1,1)$-tensor on $\gm$. Such a
 $(1,1)$-tensor is called \emph{multiplicative} \cite{crainic, PqN}. 
All it remains to show is that $J_\gm$ is completely integrable, i.e. $J_\gm$ is a multiplicative Nijenhuis tensor on $\gm$. The 
latter should presumably follow from the complete integrability of 
$J_A:TA\to TA$. This motivates the following:

\begin{quote}
\textbf{Problem 2.}
Establish a one-one correspondence between
infinitesimal multiplicative Nijenhuis tensors
on $A$ and multiplicative Nijenhuis tensors
on $\gm$.
\end{quote}

To this end, we investigate  multiplicative 
$(1, 1)$-tensors on $\Gamma$ and their infinitesimal counterparts on $A$. In particular, we show that the infinitesimal of the Nijenhuis torsion of a multiplicative $(1,1)$-tensor on $\Gamma$ is the Nijenhuis torsion of its infinitesimal counterpart on $A$. We believe that our study of multiplicative
$(1, 1)$-tensors on a Lie groupoid will be of independent interest. 

Poisson manifolds are  closely related to Lie  algebroids.
It is well known that  associated to any Poisson manifold $(X, \pi)$
there is a Lie algebroid $(T^*X)_\pi$, called the cotangent
bundle  Lie algebroid.
If the Lie algebroid $(T^*X)_\pi$ integrates to an $s$-connected
and $s$-simply connected Lie groupoid $\gm\toto X$, then
$\gm$ automatically admits a symplectic groupoid structure \cite{MX2}. In this case, the Poisson structure $(X,\pi)$ is said to be
\emph{integrable}. See \cite{CF2} (also \cite{CattaneoF})
  for the integrability conditions of  a smooth Poisson manifold.
The same situation applies to holomorphic Poisson structures
as well, and one defines
 integrable holomorphic Poisson manifolds in a similar fashion.
As an application of our general theory, we prove that 
 a holomorphic Poisson structure $(X,\pi)$,
where $\pi =\pi_R+i \pi_I$, is integrable, if and only if
$(X,\pi_R)$ and $(X,\pi_I)$ are integrable real Poisson
manifolds, recovering a theorem in
 \cite{crainic} and \cite{PqN}, which was
proved by  different methods.

The following notations are widely used in the sequel. 
For a manifold $M$, we use $q_M$ to denote the projection $TM\to M$. 
And given a complex manifold $X$, $T_\CC X$ is shorthand for the complexified tangent bundle $TX\otimes\CC$ while $T^{1,0} X$ (resp. $T^{0,1} X$) stands for the $+i$- (resp. $-i$-) eigenbundle of the almost complex structure. 
For a Lie algebroid $A$, 
the Nijenhuis torsion \cite{PN,yks} of a bundle map
$\NIJ: A\to A$ over the identity is  denoted 
$\mathcal{N}_\NIJ$, which is a section in
$\sections{\wedge^2 A^*\otimes A}$ defined by
\begin{equation}
\label{eq:nij_def}
\mathcal{N}_\NIJ(V,W)
=\lie{\NIJ V}{\NIJ W}-\NIJ(\lie{\NIJ V}{W}+\lie{V}{\NIJ W}-\NIJ \lie{V}{W}),
\qquad \forall V, W \in \sections{A} .
\end{equation}
When $A$ is the Lie algebroid $TX$ and $\NIJ: TX\to TX$ is
a $(1,1)$-tensor, the Nijenhuis torsion $\mathcal{N}_\NIJ$
is a $(2,1)$-tensor on $X$.

While the paper was in writing,  we learned   that some similar
results were also obtained independently by Ca\~{n}ez \cite{Canez}. 

\paragraph{Acknowledgments}
We would like to thank Centre \'Emile Borel and Peking
University for their hospitality while work on this project was being done.
 We also wish to thank Rui Fernandes and Alan Weinstein
for useful discussions and comments.

\section{Holomorphic Lie algebroids}

\subsection{Holomorphic Lie algebroids}

Holomorphic Lie algebroids were studied for various purposes in
the literature. See \cite{Boyom, ELW, LSX,  Hue2, Weinstein} and references cited there for details. 

%The tangent bundle $TX\to X$ of a complex manifold $X$ is naturally
%a holomorphic vector bundle. We will denote its sheaf of holomorphic sections,
%i.e. the sheaf of holomorphic vector fields, by $\Theta_X$.

%Given a holomorphic vector bundle $A\xrightarrow{p}X$,
%the sheaf of holomorphic sections $\shs{A}$ of $A\to X$ is the contravariant
%functor which associates to an open subset $U$ of $X$ the space $\shs{A}(U)$
%of holomorphic sections of $A\to X$ over $U$. Similarly, the sheaf of
% smooth sections $\shs{A}_\infty$ is the contravariant functor
% $U\to\sections{A_U}$. Clearly, $\shs{A}$ is a sheaf of $\hf{X}$-modules
%while $\shs{A}_\infty$ is a sheaf of $\cinf{X}$-modules.
%Moreover $\shs{A}$ is a subsheaf of $\shs{A}_\infty$.

By definition, a holomorphic Lie algebroid is a holomorphic vector bundle $A\to X$,
equipped with a holomorphic bundle map $A\xrightarrow{\rho}TX$,
called the anchor map, and a structure of sheaf of complex Lie algebras on $\shs{A}$,
such that
\begin{enumerate}
\item the anchor map $\rho$ induces a homomorphism of sheaves
of complex Lie algebras from $\shs{A}$ to $\Theta_X$;
\item and the Leibniz identity
\[ \lie{V}{fW}=\big(\rho(V)f\big) W+f\lie{V}{W} \]
holds for all $V,W\in\shs{A}(U)$, $f\in\hf{X}(U)$ and
  all open subsets $U$ of $X$.
\end{enumerate}
Here  $\shs{A}$ is  the sheaf of holomorphic sections  of $A\to X$
and $\Theta_X$ denotes   the sheaf of holomorphic vector fields
on $X$.

By forgetting the complex structure,
a holomorphic vector bundle $A \to X $ becomes a real (smooth) vector bundle,
and a holomorphic vector bundle map $\rho: A \to TX $
becomes a real (smooth) vector bundle map.
Assume that $A\to X$ is  a holomorphic vector bundle
 whose underlying real vector
bundle is endowed with a Lie algebroid structure $(A,\rho,\lie{\cdot}{\cdot})$
such that, for any open subset $U\subset X$,
\begin{enumerate}
\item $\lie{\shs{A}(U)}{\shs{A}(U)}\subset\shs{A}(U)$
\item and the restriction of the Lie bracket
$\lie{\cdot}{\cdot}$ to $\shs{A}(U)$ is $\CC$-linear.
\end{enumerate}
Then the restriction of $\lie{\cdot}{\cdot}$ and $\rho$
from $\shs{A}_\infty$ to $\shs{A}$ makes $A$ into
  a holomorphic Lie algebroid, where $\shs{A}_\infty$  denotes
   the sheaf of smooth sections   of $A\to X$.

The following proposition states that any holomorphic Lie algebroid
can be obtained out of such a real Lie algebroid in a unique way.

\begin{prop}[\cite{LSX}]
\label{prop:extension}
Given a structure of holomorphic Lie algebroid on the holomorphic
vector bundle $A\to X$ with anchor map $A\xrightarrow{\rho}TX$,
there exists a unique structure of real smooth Lie algebroid on
the vector bundle $A\to X$ with respect to the same anchor map
$\rho$ such that the inclusion of sheaves $\shs{A}\subset\shs{A}_\infty$
is a morphism of sheaves of real Lie algebras.
\end{prop}

By $A_R$, we   denote the underlying real Lie algebroid of a holomorphic
 Lie algebroid $A$.  
In the sequel,  by saying that a real Lie algebroid is a holomorphic Lie algebroid,
we mean that it is a holomorphic vector bundle and
its Lie bracket on smooth sections restricts to a $\CC$-linear
bracket on $\shs{A}(U)$, for all open subset $U\subset X$.

Assume that $(A\to X,\rho,\lie{\cdot}{\cdot})$ is
a holomorphic Lie algebroid. Consider
the bundle map $j: A \to A$ defining
the fiberwise complex structure on $A$. It is simple to see that
 the Nijenhuis torsion of $j$ vanishes \cite{LSX}.
Hence one can define a new (real) Lie algebroid structure on
$A$, denoted by $(A\to X,\rho_{j},\lie{\cdot}{\cdot}_j)$,
where the anchor $\rho_{j}$ is $\rho\rond j$ and
the bracket on $\sections{A}$ is
given by \cite{yks}
\[ \lie{V}{W}_{j}=\lie{jV}{W}+\lie{V}{jW}-j\lie{V}{W}, \qquad\forall V,W\in\sections{A} .\]

In the sequel, $(A\to X,\rho_{j},\lie{\cdot}{\cdot}_j)$
will be called the \emph{underlying imaginary Lie algebroid}
and denoted by $A_I$. It is known that
\begin{equation}
j: A_I\to A_R
\end{equation}
 is a Lie algebroid isomorphism \cite{yks}.

\subsection{Tangent Lie algebroids}

Recall that if $A\to M$ is a real Lie algebroid, then $TA\to TM$ is
naturally a Lie algebroid, called the tangent Lie algebroid
\cite{MX}. We recall its construction briefly below.
For details, consult \cite{MX}.

Note that for any vector bundle $E\xrightarrow{p}M$, the fibration
$TE\xrightarrow{Tp}TM$ is naturally a vector bundle. Indeed, one
has the double vector bundle:
\begin{equation}\label{tdvb}
\xymatrix{ TE \ar[r]^{Tp}\ar[d]_{q_E} & TM \ar[d]^{q_M} \\
E \ar[r]_{p} & M .}
\end{equation}

In the remainder of this section, the addition, scalar
multiplication and zero section of a fiber bundle
$E\xrightarrow{p} M$ will be denoted $+_p$, $\cdot_p$ and
$0_p$ respectively. Note that the addition, scalar
multiplication and zero section of $TE\xrightarrow{Tp}TM$ are precisely the
differentials of the addition, scalar multiplication and zero
section of $E\xrightarrow{p} M$. Hence we have
\[ \ddtz{V_t}+_{Tp}\ddtz{W_t}=\ddtz{(V_t +_p W_t)} \] if $t\mapsto V_t$
and $t\mapsto W_t$ are two paths in $E$ mapped by $p$ to the same
path in $M$ so that $\ddtz{V_t}$ and $\ddtz{W_t}$ belong to the
same fiber of $TE\xrightarrow{Tp}TM$. Similarly, for all
$\lambda\in\RR$, we have:
\[ \lambda\cdot_{Tp}\big(\ddtz{V_t}\big)=\ddtz{(\lambda\cdot_p V_t )}
.\]

%With respect to the $TP$-bundle structure, $TE\xrightarrow{Tp}TM$,
%we use $\bplus$ for addition and $\btimes$ for scalar
%multiplication. The addition $\bplus$ and the scalar
%multiplication $\bminus$ on $TE$ are precisely the tangents of the
%addition and scalar multiplication in $E$. The fiber over $v\in
%TM$ will always be denoted $(Tp)^{-1}(v)$, and the zero element of
%this fibre $T(0)(v)$. If we consider elements $\xi$ of $TE$ as
%derivatives of paths in $E$ and write
%\[ \xi = \left.\tfrac{d}{dt}V_t\right|_0 ,\]
%where $V_t$ denotes a path in $E$ defined in a neighborhood of $0\in\RR$,
% then $q_A(\xi) = V_0$ and $Tp(\xi) = \left.\tfrac{d}{dt}p(V_t)\right|_0$.
%If $\xi,\eta \in TE$ have $Tp(\xi) =Tp(\eta)$, then we can arrange that
%$\xi = \left.\tfrac{d}{dt}V_t\right|_0$ and $\eta =
%\left.\tfrac{d}{dt}W_t\right|_0$, where $p(V_t) = p(W_t)$ for all $t$ in a
%neighbourhood of $0\in\RR$. And then
%\[ \xi \dpl \eta = \left.\frac{d}{dt}(V_t + W_t)\right|_0, \qquad
%\lambda\dtimes\xi = \left.\frac{d}{dt}\lambda V_t\right|_0 .\]

Recall that, given a point $m\in M$, we have the canonical
identification \[ E_m \to T_{0_p(m)} E_m: V \mapsto \overline{V} .\]
Note that, if $V,W\in E_m$ and $\lambda\in\RR$, we have
\[ \overline{V}+_{q_E}\overline{W}=\overline{V +_p W}=
\overline{V}+_{Tp}\overline{W} \qquad \text{and} \qquad \lambda
\cdot_{q_E}\overline{V}=\overline{\lambda \cdot_p V}=
\lambda\cdot_{Tp}\overline{V} .\]

Obviously, if $V$ is a section of $E\xrightarrow{p} M$, then $TV$
is a section of $TE\xrightarrow{Tp} TM$. On the other hand, the
section $V\in\sections{E\xrightarrow{p} M}$ induces another
section $\mathcal{S}V$ of the same bundle $TE\xrightarrow{Tp} TM$
defined by the relation
\[ \mathcal{S}V(v) =T0_{p}(v)+\overline{V(m)}, \qquad\text{where }
v\in T_m M .\] Note that, for all $V,W\in\sections{E}$ and
$f\in\cinf{M}$, we have
\[ \mathcal{S}(V +_p W)=
\mathcal{S}V+_{Tp}\mathcal{S}W \qquad \text{and} \qquad
\mathcal{S}(f \cdot_p V)=(f\rond q_M)\cdot_{Tp}\mathcal{S}V .\]

%For each $m\in M$, the tangent space $T_{0_m}(E_m)$, where $0_m$
%is the zero vector of $E$ over $m$, identifies canonically with
%$E_m$; we denote the element of $T_{0_m}(E_m)$ corresponding to
%$V\in E_m$ by $\overline{V}$.
% Note that, for
%$V,W\in E_m$ and $t\in\RR$,
%\[ \overline{V} + \overline{W} = \overline{V + W} = \overline{V} \dpl \overline{W},\qquad
%t\overline{V} = \overline{tV} = t\dtimes\overline{V} .\]

%A section $V\in\sections{E}$ naturally induces a section
%$\widetilde{V}$ of $TE\xrightarrow{Tp}TM$ by
%\[ \widetilde{V}(v) = T(0)(v)+ \overline{V(m)}  , \qquad \forall v\in T_m M .\]
%Note that
%\[ \widetilde{(V + W)} = \widetilde{V}\dpl\widetilde{W},\qquad
%\widetilde{(fV)} = (f\rond q)\dtimes\widetilde{V} ,\] for
%$V,W\in\sect, f\in C^{\infty}(M)$.

%When $E$ is a trivialized vector bundle $M\times F\lon M$, where
%$F$ is a vector space, we denote elements of $TE = TM\times F\times F$ by
%$\xi = (v,V,V')$,
%where if $m_t$ is a path in $M$ with $v = \ddt{m_t}$, we identify $(v,V,V')$
%with $\ddt{(m_t,V + tV')}$.
%The vector bundle structure of $Tp: TM\times
%F\times F\lon TM: (v,V,V') \mapsto v$ is then given by
%\begin{align}
%(v,V,V')\dpl (v,W,W') & =  (v,V+W,V'+W'), \nonumber \\
%t\dtimes(v,V,V') & = (x,tV, tV'), \label{eq:trivTA} \\
%T(0)(v) & =  (v,0,0), \nonumber
%\end{align}
%where $v\in TM,\ V, V', W, W'\in F,\ t\in\RR$.

If $E$ is a trivialized vector bundle $M\times F\lon M$,
a section is equivalent to a map $V: M\to F$.
Then $\mathcal{S}V$ is the section of $Tp: TM\times
F\times F\lon TM$ given by
\[ TM\to TM\times F\times F: v_m\mapsto (v_m, 0, V(m)) .\]

On the other hand, if $V\in \sections{E}$, $TV$ defines  another section
of $Tp: TE\to TM$. Under the trivialization $p: M\times F\to M$,
$TV$ is the section
\[ TM\to TM\times F\times F: v_m\mapsto (v_m, V(m), v_m (V)) .\]

It is easy to see that $\sections{TE\xrightarrow{Tp}TM}$ is
generated by the sections of type $TV$ and $\mathcal{S}V$.

Taking $E=TM$ and $p=q_M$ in Diagram~\eqref{tdvb}, we obtain
the double vector bundle
\[ \xymatrix{ T(TM) \ar[r]^{Tq_M}\ar[d]_{q_{TM}} & TM \ar[d]^{q_M} \\
TM \ar[r]_{q_M} & M .} \]
As is clearly shown by this last diagram, $T(TM)$ is endowed with
two different vector bundle structures over the same base $TM$.
However, there exists a canonical isomorphism of double vector
bundles \[ \sigma : T(TM)\to T(TM), \]
called the \emph{flip map}, which exchanges these two
vector bundle structures $T(TM)\to TM$:
\[ \xymatrix{ T(TM) \ar[rrd]^{q_{TM}} \ar[ddr]_{Tq_{M}}
\ar[dr]_{\sigma_M} & & \\
& T(TM) \ar[r]_{Tq_M}\ar[d]^{q_{TM}} & TM \ar[d]^{q_M} \\
& TM \ar[r]_{q_M} & M .} \]

The following theorem was proved in \cite{MX}.

\begin{thm}\label{thm:TA}
Let $(A\xrightarrow{p}M,\rho,[\cdot,\cdot])$ be a Lie algebroid.
Then there exists a (unique) Lie algebroid structure on the vector
bundle $TA\xrightarrow{Tp}TM$, with anchor $\rho_T:=\sigma_\rond\rho_*$,
whose bracket on sections satisfies the relations
\begin{equation}
\label{eq:TS} \lie{TV}{TW}=T\lie{V}{W}
\qquad
\lie{TV}{\mathcal{S}W}=\mathcal{S}\lie{V}{W}
\qquad \lie{\mathcal{S}V}{\mathcal{S}W}=0 
\end{equation}
 for any $V,W\in\sections{A\xrightarrow{p}M}$.
\end{thm}

\subsection{Equivalent definition of holomorphic Lie algebroids}

\begin{prop}\label{prop:tgtalg}
Let $(A,\rho, [\cdot,\cdot] )$ be a real Lie algebroid,
where $A \to X $ is a holomorphic vector bundle.
 The following are equivalent:
\begin{enumerate}
\item $(A,\rho, [\cdot,\cdot] )$ is a holomorphic Lie algebroid;
\item the map
\begin{equation}
\label{eq:TA}
 \xymatrix{ TA \ar[d] \ar[r]^{J_A} & TA \ar[d] \\ TX \ar[r]_{J_X} & TX } 
\end{equation}
defines a Lie algebroid isomorphism,
where $J_A$ and $J_X $ denote the almost complex structures on $A$
and $X$ respectively.
\end{enumerate}
\end{prop}

\begin{proof}
First of all, it is simple to check, using  the local
trivialization of the holomorphic bundle $A\to X$, that
$(J_A, J_X)$ is indeed a bundle map.
It follows from Theorem~\ref{thm:TA} and Eq.~\eqref{eq:TS} that
Diagram~\eqref{eq:TA} is a Lie algebroid isomorphism if, and only if,
\begin{equation} \label{eq:pf3}
\rho_{T} \rond J_A =  J_{X*} \rond \rho_T
\end{equation}
and
\begin{equation} \label{eq:pf4}
\left\{ \begin{aligned} 
& \lie{ J_A \rond TV \rond J_X\inv }{ J_A \rond TW \rond J_X\inv } 
= J_A \rond [TV,TW] \rond J_X\inv \\
& \lie{ J_A \rond TV \rond J_X\inv }{ J_A \rond \widetild{W} \rond J_{X}\inv } 
= J_A \rond \widetild{\lie{V}{W}} \rond J_X^{-1} \\ 
& \lie{ J_A \rond \widetild{V} \rond J_X\inv }{ J_A \rond \widetild{W} \rond
J_X\inv } = 0 
\end{aligned} \right. 
\end{equation}
for any open subset $U\subseteq X$ and all $V,W\in\shs{A}(U)$.

The proposition now becomes an immediate consequence of the following
two observations.

\textsl{(1) The anchor map $\rho: A \to TX $ is a holomorphic map
if, and only if, Eq.~\eqref{eq:pf3} holds.}

The anchor map $\rho: A \to TX $ is a holomorphic map if, and only if,
\begin{equation} \label{eq:pf2}
\xymatrix{ TA\ar[r]^{J_A} \ar[d]_{\rho_*} & TA \ar[d]^{\rho_*} \\ T(TX)
\ar[r]_{J_{TX}} & T(TX) }
\end{equation}
commutes, where $J_{TX}$ is the almost complex structure
on $TX$ induced from the one on $X$.
It is known that $J_{TX}$ is given by the following equation:
\[ J_{TX} = \sigma \rond  J_{X*}  \rond \sigma : T(TX)\to T(TX) .\]
Conjugating the second horizontal line in Eq.~\eqref{eq:pf2} by the
flip map $\sigma_X $, we  see that Diagram~\eqref{eq:pf2} commutes if, and only if,
\[ \xymatrix{ TA\ar[r]^{J_A} \ar[d]_{\rho_{T}} & TA \ar[d]^{\rho_{T}} \\ 
T(TX) \ar[r]_{J_{X*}} & T(TX) } \]
commutes, and the latter is precisely Eq.~\eqref{eq:pf3}.

\textsl{(2) The Lie bracket on $\sections{A_U}$ restricts to a $\CC$-linear bracket on $\shs{A}(U)$
if, and only if, Eq.~\eqref{eq:pf4} is satisfied for all $V,W\in\shs{A}(U)$.}

A smooth map between holomorphic manifolds is holomorphic if, and
only if, its differential intertwines the corresponding almost complex structures.
In particular, a smooth section $V\in \sections{A_U}$ is holomorphic if, and only if,
\begin{equation}\label{eq:hp6}
J_A \rond T V \rond J_X^{-1} = TV .
\end{equation}

For all $V,W \in \shs{A} (U)$, we have, on the one hand,
\[ J_A \rond T [V,W] \rond J_X\inv = J_A \rond  [TV,TW]\rond J_X\inv \]
and, on the other hand according to Eq.~\eqref{eq:hp6},
\[ T[V,W]=[TV,TW]=[J_A\rond TV\rond J_X\inv,J_A\rond TW\rond J_X\inv] .\]
Hence the bracket $[V,W]$ is holomorphic if, and only if,
\[ [ J_A \rond TV \rond J_X\inv, J_A \rond TW \rond
J_{X}\inv] = J_A \rond [TV,TW] \rond J_X\inv \]
holds for all $V,W \in \shs{A}(U)$.

Secondly it is simple to check that
\begin{equation}\label{eq:hp7}
J_A \rond \widetild{V} \rond J_X^{-1} = \widetild{(jV)} 
.\end{equation}

For any $V,W \in \shs{A}(U)$, we have, on the one hand,
\begin{align*}
\widetild{\lie{V}{jW}} &= \lie{TV}{\widetild{(jW)}}  && \text{(by Eq.~\eqref{eq:TS})} \\
&= [J_A\rond TV\rond J_X\inv,J_A\rond\widetild{W}\rond J_X\inv]
&& \text{(by Eqs.~\eqref{eq:hp6}~and~\eqref{eq:hp7})}
\end{align*} 
and, on the other hand according to Eq.~\eqref{eq:hp7},
\[ \widetild{(j\lie{V}{W})} = J_A \rond \widetild{\lie{V}{W}} \rond J_X\inv .\]
It thus follows that ${\lie{V}{jW}} = {j\lie{V}{W}}$ if, and only if,
\[ [ J_A \rond TV \rond J_X\inv, J_A \rond \widetild{W} \rond
J_{X}\inv] = J_A\rond \widetild{\lie{V}{W}}\rond J_X\inv \] 
holds for all $V,W\in\shs{A}(U)$.

Finally, according to Eq.~\eqref{eq:hp6}, the relation
$\lie{ J_A \rond \widetild{V} \rond J_X\inv }
{ J_A \rond \widetild{W} \rond J_X\inv}=0$
holds for any $V,W\in\shs{A}(U)$.
\end{proof}

\subsection{Holomorphic Poisson manifolds}
% and their associated holomorphic Lie groupoids}

\begin{defn}
A holomorphic Poisson manifold is a  complex manifold $X$ equipped
with a holomorphic bivector field $\pi$ 
(i.e. $\pi\in \gm (\wedge^2 T^{1, 0}X)$
such that $\delbar\pi=0 $), satisfying the equation $\schouten{\pi}{\pi}=0$.
\end{defn}

Since $\wedge^2 T_{\CC}X = \wedge^2 TX \oplus i \wedge^2 TX$,
for any $\pi \in \sections{\wedge^2 T_{\CC}X}$, we can write
$\pi=\pire+i\piim$, where $\pire$ and $\piim\in\sections{\wedge^2 TX}$
are (real) bivector fields on $X$. The following
result was proved in \cite{LSX}.

%Both $\pire$ and $\piim$ define
%brackets $\pbre{\cdot}{\cdot}$ and $\pbim{\cdot}{\cdot}$
%on $\cinf{M,\RR}$ in the standard way.
%These extend to $\cinf{M,\CC}$ by $\CC$-linearity.
%The next lemma describes such an extension.

\begin{thm}[\cite{LSX}]
\label{PNGC} 
Given a complex manifold $X$ with associated
almost complex structure $J$, the following are equivalent:
\begin{enumerate}
\item\label{pq1} $\pi=\pire+i\piim\in\sections{\wedge^2T^{1,0} X}$
is a holomorphic Poisson bivector field;
\item\label{pq2} $(\piim,J)$ is a Poisson Nijenhuis structure on $X$ and
$\pire\diese=\piim\diese\rond J^*$;
\end{enumerate}
\end{thm}

As a consequence, both $\pire$ and $\piim$ are Poisson tensors
and they constitute a biHamiltonian system \cite{PN}.
Both $\pire$ and $\piim$ define
brackets $\pbre{\cdot}{\cdot}$ and $\pbim{\cdot}{\cdot}$
on $\cinf{M,\RR}$ in the standard way.
These extend to $\cinf{M,\CC}$ by $\CC$-linearity.
 In a local chart $(z_1=x_1+i y_1,\cdots,z_n=x_n+i y_n)$ of complex coordinates of $X$,
we indeed have
\begin{align*}
\pbre{x_i}{x_j}&=\tfrac{1}{4}\Re\pb{z_i}{z_j}, & \pbim{x_i}{x_j}&=\tfrac{1}{4}\Im\pb{z_i}{z_j}, \\
\pbre{y_i}{y_j}&=-\tfrac{1}{4}\Re\pb{z_i}{z_j}, & \pbim{y_i}{y_j}&=-\tfrac{1}{4}\Im\pb{z_i}{z_j}, \\
\pbre{x_i}{y_j}&=\tfrac{1}{4}\Im\pb{z_i}{z_j}, & \pbim{x_i}{y_j}&=-\tfrac{1}{4}\Re\pb{z_i}{z_j}.
\end{align*}
Here $\Re$ and $\Im$ denote the real and imaginary parts of a complex number.

Now   we consider the cotangent  bundle Lie algebroid
of  a holomorphic Poisson manifold and identify its associated
real and imaginary  Lie algebroids.
Assume that $(X, \pi)$ is a holomorphic Poisson manifold, where
$\pi=\pi_R+i \pi_I\in \gm (\wedge^2 T^{1, 0}X)$. Let
$A=(T^*X)_\pi$ be its corresponding holomorphic Lie algebroid,
which can be  defined in a similar way as in the smooth case.
To be more precise, let $\Phi$ and $\Psi$, respectively, be the
holomorphic bundle maps
\[ \Phi: TX \to T^{1, 0}X, \ \Phi=  \frac{1}{2} (1-iJ)  \]
and
\[ \Psi: T^*X\to (T^{1, 0}X)^*, \Psi=1-iJ^* ,\]
where $J$ is the almost complex structure on $X$.
Define the anchor $\rho: (T^*X)_\pi\to TX$ to be
$\rho=\Phi^{-1}\rond \pi^{\#} \rond \Psi$ and the bracket
\[ [\alpha, \beta]_\pi =  L_{\rho \alpha } \beta - L_{\rho \beta }
   \alpha  -   {\rm d} (\rho \alpha,   \beta ) \]
 $\forall \alpha,\beta  \in \Gamma ( T^*X|_U)$ holomorphic.
One easily sees that $(T^* X)_{\pi} $ is a  holomorphic
Lie algebroid. Next proposition describes its real and
imaginary Lie algebroids.

\begin{prop}[\cite{LSX}]
\label{prop:3.7}
Let $(X, \pi )$ be a holomorphic Poisson manifold,
 where
$\pi=\pi_R+i \pi_I\in \gm (\wedge^2 T^{1, 0}X)$. Then
the underlying real and
imaginary Lie algebroids of $(T^*X)_\pi$ are isomorphic
to $(T^*X)_{4\pi_R}$ and $(T^*X)_{4\pi_I}$, respectively.
\end{prop}

\section{Holomorphic Lie groupoids and holomorphic symplectic groupoids}

\subsection{Multiplicative (1,1)-tensors on Lie groupoids}

Recall that a skew-symmetric $(k,1)$-tensor on a smooth manifold
$M$ can be seen either as a section of
the vector bundle $\wedge^{k} T^*M \otimes TM \to M$,
or as a bundle map $\wedge^{k} TM\to TM$
over the identity on $M$.
If $\Gamma \toto M$ is a Lie groupoid, then $T\Gamma \toto TM$ is
naturally a Lie groupoid, called the tangent Lie groupoid \cite{MX}.
Indeed for any $k\geq 1$, this Lie groupoid structure extends naturally to a Lie groupoid $\wedge^k T\Gamma \toto \wedge^kT M$, whose 
source, target, inverse and product maps are  given by
$\wedge^k T s$, $\wedge^k T t$, $\wedge^k T\inverse$ and
$\wedge^k T\product$, respectively.
Here $s$, $t$, $\inverse$ and $\product$
are the source, target, inverse and product maps of the
groupoid $\Gamma\toto M$.

\begin{defn}
A multiplicative $(k,1)$-tensor $\NIJ$ on a Lie groupoid
$\Gamma\toto M$
is a pair $(\NIJ_\Gamma,\NIJ_M)$ of $(k,1)$-tensors on
 $\Gamma$ and $M$ respectively
such that
\[ \xymatrix{
\wedge^k T\Gamma \dar[d] \ar[r]^{\NIJ_\Gamma} & T\Gamma \dar[d] \\
\wedge^k TM \ar[r]_{\NIJ_M} & TM
} \]
is a Lie groupoid morphism.
\end{defn}

\begin{rmk}
It is simple to see that if 
$(\NIJ_\Gamma,\NIJ_M)$ is  a multiplicative $(k,1)$-tensor,
$\NIJ_M$ is completely determined by $\NIJ_\Gamma$, which
is the restriction of $\NIJ_\Gamma$ to the unit space $M$.
Hence we often use $\NIJ_\Gamma$
 to denote a multiplicative $(k,1)$-tensor. 
\end{rmk}

%Recall that the Nijenhuis tensor $\mathcal{N}_\NIJ$ of a $(1,1)$-tensor
%$\NIJ$ on a manifold $M$ is the $(2,1)$-tensor on $M$ defined by
%\beq{eq:nij_def} \mathcal{N}_\NIJ(X,Y) = [\NIJ X,\NIJ Y] - \NIJ([\NIJ X,Y]+
%[X,\NIJ Y]-\NIJ[X,Y]) \eeq
%for all $X,Y\in\vf(M)$.

\begin{prop}\label{prop:mult_nij}
If $(\NIJ_\Gamma,\NIJ_M)$ is a multiplicative $(1,1)$-tensor on a
Lie groupoid $\Gamma\toto M$, then $(\mathcal{N}_{\NIJ_\Gamma},\mathcal{N}_{\NIJ_M})$
is a multiplicative $(2,1)$-tensor on $\Gamma$.
\end{prop}

The proof requires two lemmas.
The first one, a general fact about Nijenhuis tensors,
is a straightforward consequence of Eq.~\eqref{eq:nij_def}.
A $(k,1)$-tensor $\NIJ$ on a manifold $N$ is said to be \emph{tangent
to a submanifold $S\subset N$} if $\NIJ$ maps $\wedge^k TS$ to $TS$.
Any $(k,1)$-tensor $\NIJ$ tangent to $S$
induces by restriction a $(k,1)$-tensor on the submanifold $S$.

The following lemma can be easily verified.

\begin{lem}\label{lem:nij_tgt}
If $S\subset N$ is a submanifold, and $\NIJ$ is $(1,1)$-tensor
tangent to $S$, then $\mathcal{N}_\NIJ$ is tangent to $S$.
Moreover the restriction of $\mathcal{N}_\NIJ$ to $S$
is the Nijenhuis tensor of the restriction
of $\NIJ$ to $S$.
\end{lem}

The second lemma is a general fact regarding Lie groupoids,
the proof of which is left to the reader.
Recall that the \emph{graph of the multiplication}
of a given Lie groupoid $\Gamma\toto M $ is the 
submanifold $\Lambda_\Gamma$
\[ \Lambda_\Gamma = \genrel{(\gamma_1,\gamma_2,\gamma_1\gamma_2)}
{\forall \gamma_1,\gamma_2\in\Gamma\text{ s.t. }t(\gamma_1)= s(\gamma_2)} \]
of $\Gamma\times\Gamma\times\Gamma$.
Given a map $\Psi:R\to R'$, we denote by $\Psi^{(3)}$ the map
\[ \Psi^{(3)}:R\times R\times R\to R'\times R'\times R':(x,y,z)\mapsto\big(\Psi(x),\Psi(y),\Psi(z)\big) .\]

\begin{lem}\label{lem:groids_graph}
\begin{enumerate}
\item\label{65a} The graphs of the multiplications of the Lie groupoids
$\Gamma\toto M$, $T\Gamma\toto TM$ and $\wedge^2 T\Gamma\toto\wedge^2 TM$
are related as follows
\[ \Lambda_{T\Gamma}=T\Lambda_{\Gamma} \qquad\qquad \Lambda_{\wedge^2 T\Gamma}=\wedge^2 T\Lambda_{\Gamma} .\]
\item\label{65b} A smooth map $\Psi$ from a Lie groupoid
$\Gamma\toto M$ to another Lie groupoid $\Gamma'\toto M'$
is a groupoid homomorphism if, and only if,
$\Psi^{(3)}$ maps $\Lambda_{\Gamma}$ to $\Lambda_{\Gamma'}$.
\end{enumerate}
\end{lem}

\begin{proof}[Proof of Proposition~\ref{prop:mult_nij}]
By Lemma~\ref{lem:groids_graph}(\ref{65b}), since the pair
$(\NIJ_\Gamma,\NIJ_M)$ is a Lie groupoid homomorphism,
$\NIJ_\Gamma^{(3)}:T\Gamma^{3}\to T\Gamma^{3}$ maps
$\Lambda_{T\Gamma}$ to itself.
By Lemma~\ref{lem:groids_graph}(\ref{65a}), the latter
is isomorphic to $T\Lambda_\gm$.
Hence $\NIJ_\Gamma^{(3)}$ is tangent to $\Lambda_\Gamma$.

By Lemma~\ref{lem:nij_tgt}, it follows that $\mathcal{N}_{\NIJ^{(3)}_\Gamma}$
is also tangent to $\Lambda_{\Gamma}$.
A simple computation yields that
\[ \mathcal{N}_{\NIJ_\Gamma^{(3)}}=\mathcal{N}_{\NIJ_\Gamma}^{(3)} .\]

Therefore $\mathcal{N}_{\NIJ_\Gamma}^{(3)}$ maps
$\wedge^2 T\Lambda_\Gamma$ to $T\Lambda_\Gamma$.
According to Lemma~\ref{lem:groids_graph} again,
this amounts to saying that $(\mathcal{N}_{\NIJ_\Gamma},\mathcal{N}_{\NIJ_M})$
is indeed a Lie groupoid homomorphism.
\end{proof}

\subsection{Nijenhuis torsion of multiplicative (1,1)-tensors}

Recall that, for any vector bundle $E\xrightarrow{p}M$, one has
the double vector bundle \eqref{tdvb}. Similarly, one has
the double vector bundle
\begin{equation}\label{ttdvb}
\xymatrix{ \wedge^2 TE \ar[r]^{Tp}\ar[d]_{q_E} & \wedge^2 TM \ar[d]^{q_M} \\
E \ar[r]_{p} & M .}
\end{equation}

In particular, when $E$ is the tangent bundle $TM\xrightarrow{p}M$, 
we have the double bundles
\begin{equation}\label{ttm}
\xymatrix{ T(TM) \ar[r]^{Tp}\ar[d]_{{q_{TM}}} & TM \ar[d]^{q_M} \\
TM \ar[r]_{p} & M }
\end{equation}
and
\begin{equation}\label{tttm}
\xymatrix{ \wedge^2 T(TM) \ar[r]^{Tp}\ar[d]_{q_{TM}} & \wedge^2 TM \ar[d]^{q_M} \\
TM \ar[r]_{p} & M .}
\end{equation}

There is an extension of the canonical flip map $\sigma: T(TM)\to T(TM)$,
denoted by $\sigma^{(1)}$ in this section, to
$T(\wedge^2 TM)$, which is defined as follows.
Let $\mu\in T_{e_1 \wedge e_2}(\wedge^2 TM)$ be any
tangent vector, where $e_1, e_2\in T_m M$. Write
\[ \mu= \left.\tfrac{d}{dt}e_1 (t)\wedge e_2 (t) \right|_0 ,\]
where $e_1(t),e_2(t)\in T_{m(t)}M$.
Define
\[ \sigma^{(2)}(\mu)= \left. \sigma^{(1)} \tfrac{d}{dt} e_1(t)
\right|_0 \wedge \left. \sigma^{(1)} \tfrac{d}{dt} e_2 (t) \right|_0 .\]
Then $\sigma^{(2)}(\mu)$ is a vector in $\wedge^2 T_v(TM)$.
Hence $\sigma^{(2)}$ maps $T_{e_1 \wedge e_2}(\wedge^2 TM)$
to $\wedge^2 T_v (TM)$, where $v=(Tp)(e_1 \wedge e_2)\in T_mM$.
Note that $v$ changes according to $e_1 \wedge e_2$.

\begin{lem}
\label{lem:6.5}
\begin{enumerate}
\item $\sigma^{(1)}$ is an isomorphism of double vector bundle
\eqref{ttm}, which interchanges the horizontal and vertical
bundle structures and induces the identities on the side bundles.
\item $\sigma^{(2)}$ is an isomorphism from the double vector bundle
\begin{equation*}
\xymatrix{ T(\wedge^2 TM) \ar[r]^{q_{\wedge^2 TM}}\ar[d]_{Tp} &
 \wedge^2 TM \ar[d]^{q_M} \\
TM \ar[r]_{p} & M }
\end{equation*}
to the double vector bundle \eqref{tttm}, 
which induces the identities on the side bundles.
\end{enumerate}
\end{lem}

Let $\NIJ: TM\to TM$ be a $(1,1)$-tensor on $M$. Then $T\NIJ$
is a morphism of  horizontal  vector bundle $T(TM)\stackrel{Tp}{\to} TM$ 
to itself over the identity.
Define
\[ \mathbb{T}\NIJ = \sigma^{(1)} \rond T\NIJ \rond (\sigma^{(1)})^{-1} .\]
From Lemma~\ref{lem:6.5},
 it follows that $\mathbb{T}\phi$ is a
morphism of vertical  vector bundle $q_{TM}:T(TM)\to TM$ to itself over
the identity.

Similarly, if $\NIJ:\wedge^2TM\to TM$ is  a $(2,1)$-tensor on $M$,
then \[ \mathbb{T}\phi = \sigma^{(1)} \rond T\NIJ \rond (\sigma^{(2)})^{-1} \]
is a morphism from the vector bundle $q_{TM}:\wedge^2 T(TM)\to TM$
to the vector bundle $q_{TM}:T(TM)\to TM$ over the identity.

\begin{prop}\label{pro:difftensors}
If $\NIJ$ is a $(1,1)$-tensor (resp. $(2,1)$-tensor) on a manifold $M$, then
$\mathbb{T}\NIJ$ is a $(1,1)$-tensor (resp. $(2,1)$-tensor) on the manifold $TM$.
Moreover, we have
\begin{equation}\label{eq:pre_nij}
\mathcal{N}_{\mathbb{T}\NIJ} = \mathbb{T}\mathcal{N}_{\NIJ}
.\end{equation}
\end{prop}

\begin{proof}
It remains to prove Eq.~\eqref{eq:pre_nij}.
Note that for any section $V$ of $p:TM\to M$, $TV$ is a section of
$Tp:T(TM)\to TM$, which is a tangent Lie algebroid.
Let \[ \mathbb{T}V=\sigma^{(1)}\rond TV .\]
Then $\mathbb{T}V$ is a section of the bundle $q_{TM}:T(TM)\to TM$, i.e.
a vector field on $TM$. Considering $p:TM\to M$ as a Lie algebroid
and $Tp: T(TM) \to TM$ as its tangent Lie algebroid,
according to Eq.~\eqref{eq:TS}, we have $[TV,TW]=T[V,W]$, 
where the left hand side of
the equation refers to the Lie bracket of the tangent Lie algebroid
$T(TM)\to TM$. Hence
\begin{equation}
\label{eq:T1}
 \mathbb{T}[V,W]=\sigma^{(1)}\rond T[V,W]=\sigma^{(1)}\rond
[TV,TW]=[\mathbb{T}V,\mathbb{T}W]. 
\end{equation} 
Here the  bracket $[\mathbb{T}V,\mathbb{T}W]$ on the right hand
side  refers to the bracket
on the vector fields $\vf(TM)$, and the last equality follows from
the definition of tangent Lie algebroid.

Now if $\NIJ$ is a $(1,1)$-tensor on $M$, we have
\begin{equation}\label{eq:hatT}
(\mathbb{T}\NIJ)(\mathbb{T}V) = \sigma^{(1)}\rond T\NIJ \rond
(\sigma^{(1)})^{-1} \rond \sigma^{(1)} \rond TV = \sigma^{(1)}
 \rond T\NIJ \rond TV =\sigma^{(1)} \rond T\NIJ(V) = \mathbb{T}(\NIJ V)
.\end{equation}
Similarly, for any $(2,1)$-tensor $\psi$ on $M$ and any sections
$V,W$ of $TM\to M$, we have
\begin{equation}\label{eq:hatT2}
(\mathbb{T}\psi)(\mathbb{T}V,\mathbb{T}W) = \mathbb{T}\big(\psi(V,W)\big)
.\end{equation}

Using Eqs.~\eqref{eq:T1} and \eqref{eq:hatT2}, a simple computation leads to 
\[ {\mathcal N}_{\mathbb{T} \NIJ} ( \mathbb{T}V , \mathbb{T}W )
= (\mathbb{T} {\mathcal N}_\NIJ) ( \mathbb{T}V,\mathbb{T}W ) .\]

%\[
%\begin{array}{rcll}
%{\mathcal N}_{\mathbb{T} \NIJ} ( \mathbb{T}V , \mathbb{T}W )&=& (\mathbb{T} \NIJ)^2 ( [ \mathbb{T}V , \mathbb{T}W ] )
%  +  [ (\mathbb{T} \NIJ) \rond \mathbb{T}V ,(\mathbb{T} \NIJ) \rond \mathbb{T}W ] \\ & &  -
%  \mathbb{T}\NIJ  [ (\mathbb{T} \NIJ) \rond \mathbb{T}V , \mathbb{T}W ]  -  \mathbb{T}\NIJ [ \mathbb{T}V , (\mathbb{T} \NIJ ) \rond \mathbb{T}W ] \big) & \\
%       &=& (\mathbb{T} \NIJ)^2 ( [ \mathbb{T}V , \mathbb{T}W ] )  + [ \mathbb{T} (\NIJ (V)) ,\mathbb{T} (\NIJ W) ]  \\ & &
%         - (\mathbb{T}\NIJ)  \big( [ (\mathbb{T} \NIJ(V) , \mathbb{T}W ]  - [  \mathbb{T}V , \mathbb{T} (\NIJ W) ] \big)& \mbox{by Eq. (\ref{eq:hatT})} \\
%        &=&  (\mathbb{T} \NIJ)^2 ( \mathbb{T}[ V , W ] )  + \mathbb{T}[  \NIJ (V) , \NIJ W ]  & \\ & &
%         - (\mathbb{T} \mathbb{\NIJ} )  \mathbb{T}[  \NIJ(V) , W ])  - (\mathbb{T} \mathbb{\NIJ} ) \mathbb{T} [ V , \NIJ W ] )&  \mbox{ by Lemma \ref{lem:hat} (1)}  \\
%         &=&  \mathbb{T} (\NIJ^2 ( [ V , W ] )  + \mathbb{T}[  \NIJ V , \NIJ W ]  & \\ & &
%         - \mathbb{T} (\NIJ [  \NIJ(V) , W ])  - \mathbb{T}( \NIJ  [ V , \NIJ W ] )&   \mbox{by Eq. (\ref{eq:hatT})} \\
%         &=&  \mathbb{T} \big( {\mathcal N}_{ \NIJ}(V,W ) \big) \\
%         &= &   (\mathbb{T} {\mathcal N}_\NIJ) ( \mathbb{T}V,\mathbb{T}W  )   & \mbox{ by Eq. (\ref{eq:hatT2})} \end{array}
%\]

On the other hand, using local coordinates, it is easy to see that
for any $u\in TM$  with $u\neq 0$ and $v\in T_u(TM)$,
there always exists a section $V$ of $TM\to M$ such that $\mathbb{T}V$ passes
through $v$ at $u$.
Since both $\mathcal{N}_{\mathbb{T}\NIJ}$ and $\mathbb{T}\mathcal{N}_\NIJ$
are $(2,1)$-tensors on $TM$, it follows that they coincide
at all points of $TM$ except for the zero section of
$TM$. Hence they must be equal at all points by continuity.
\end{proof}

%\subsection{Main theorem}

Now let $\gm\toto M$ be a Lie groupoid with Lie algebroid $A$.
By definition $A$ is identified with the subbundle $T^s_M\gm$
of $T\gm|_M$. In the sequel, we use this identification implicitly.

Let $(\NIJ_\gm,\NIJ_M)$ be a multiplicative $(1,1)$-tensor on
$\gm$. Then by Proposition~\ref{pro:difftensors},
 $\mathbb{T}\NIJ_\gm$ is a $(1,1)$-tensor
on the manifold $T\gm$. Similarly if $(\psi_\gm,\psi_M)$ is a
multiplicative $(2,1)$-tensor on $\gm$, then $\mathbb{T}\psi_\gm$ is a
$(2,1)$-tensor on the manifold $T\gm$.

\begin{lem}
Let $(\NIJ_\gm,\NIJ_M)$ and $(\psi_\gm,\psi_M)$
be a multiplicative $(1,1)$-tensor and $(2,1)$-tensor
on a Lie groupoid $\gm\toto M$, respectively.
Then both $\mathbb{T}\NIJ_\gm$ and $\mathbb{T}\psi_\gm$
are tangent to the submanifold $T^s_M \gm\subset T\gm$.
Hence they define a $(1,1)$-tensor and $(2,1)$-tensor
on the submanifold $T^s_M \gm$ of $T\gm$.
\end{lem}

\begin{proof}
By definition,
$\mathbb{T}\NIJ_\gm=\sigma^{(1)}\rond T\NIJ\rond 
(\sigma^{(1)})^{-1}$.
It is well-known that 
$\sigma^{(1)}:  T_{TM}^{Ts}T\gm\to T (T_M^{s}\gm )$ is an isomorphism.
 Since $\NIJ: T\gm\to T\gm$ is a
groupoid morphism of the tangent groupoid $T\gm\toto TM$, it follows that
$T_{TM}^{Ts}T\gm$ is stable under $T\NIJ$.
Hence it follows that $T (T_M^{s}\gm )$ is stable under
$\mathbb{T}\NIJ$. That is, $\mathbb{T}\NIJ_\gm$
is tangent to the submanifold $T^s_M \gm\subset T\gm$.

Similarly, one proves that $\mathbb{T}\psi_\gm$ is also tangent
to $T^s_M \gm$.
\end{proof}

Now we introduce

\begin{defn}
If $(\NIJ_\gm, \NIJ_M)$ (resp. $(\psi_\gm, \psi_M)$) is
a multiplicative $(1,1)$-tensor (resp. $(2,1)$-tensor) on
$\gm$, we define $\LLie \NIJ_\gm$ (resp. $\LLie \psi_\gm$)
to be the restriction of $\mathbb{T}\NIJ_\gm$ (resp. $\mathbb{T}\psi_\gm$)
to $A$ (being identified with $T^s_M \gm$ and
considered as a submanifold of $T\gm$).
\end{defn}

The main result of this section is the following:

\begin{thm}\label{th:Lie_deriv}
Let $\gm\toto M$ be a Lie groupoid with Lie algebroid $A$.
If $\NIJ_\gm$ is a multiplicative $(1,1)$-tensor on $\Gamma$,
then ${\LLie (\NIJ_\gm)}$ is a $(1,1)$-tensor on $A$ (as a manifold),
and ${\LLie(\mathcal{N}_{\NIJ_\gm})}$
is a $(2,1)$-tensor on $A$ (as a manifold). Moreover we have
\begin{equation}\label{eq:theo_nij}
\mathcal{N}_{ {\LLie (\NIJ_\gm)} } =  {\LLie(\mathcal{N}_{\NIJ_\gm})}
.\end{equation}
\end{thm}

\begin{proof}
In Lemma~\ref{lem:nij_tgt}, taking $N=T\gm$ and $S=T_M^s\gm$, we obtain
that $\mathcal{N}_{\LLie \NIJ_\gm }=\mathcal{N}_{\mathbb{T}\NIJ_\gm}|_{T_M^s \gm}$.
The latter is equal to 
$\mathbb{T}{\mathcal N}_{\NIJ_\gm}|_{T_M^s \gm}$ according to 
Proposition~\ref{pro:difftensors},
which is $\LLie(\mathcal{N}_{\NIJ_\gm})$ by definition.
\end{proof}

\subsection{Infinitesimal multiplicative (1,1)-tensors}

\begin{defn}
Let $(A,\rho,[\cdot,\cdot])$ be a real Lie algebroid.
An infinitesimal multiplicative $(1,1)$-tensor on $A$
is a pair $(\NIJ_A,\NIJ_M)$ of $(1,1)$-tensors on $A$ and $M$
such that the following diagram  
\begin{equation}\label{eq:Axxx}
\xymatrix{ TA \ar[r]^{\NIJ_A} \ar[d] & TA \ar[d] \\ TM \ar[r]_{\NIJ_M} & TM } 
\end{equation}
is a Lie algebroid homomorphism, where $TA\to TM$ is the tangent
Lie algebroid of $A\to M$.
\end{defn}

\begin{rmk}
It is simple to see that $\NIJ_M$ is completely determined by 
$\NIJ_A$, which is equal to the restriction 
of $\NIJ_A$ to the zero section $M\subset A$. In the sequel,
we sometimes simply use $\NIJ_A$ to denote an 
infinitesimal multiplicative $(1,1)$-tensor on $A$.
\end{rmk}

In particular, for any holomorphic Lie algebroid $A\to X$, 
the almost complex structures $(J_A,J_X)$ is an infinitesimal 
multiplicative $(1,1)$-tensor on $A$ according to Proposition~\ref{prop:tgtalg}.

Let $(\NIJ_{\Gamma},\NIJ_M)$ be a multiplicative $(1,1)$-tensor on a Lie groupoid $\Gamma\toto M$.
According to Theorem~\ref{th:Lie_deriv}, $\LLie(\NIJ_\Gamma)$ is a $(1,1)$-tensor on $A$.
Moreover, by construction, it is clear that $\LLie(\NIJ_\Gamma)$ is a homomorphism of
the tangent Lie algebroid $TA\to TM$ to itself. Hence the assignment 
\begin{equation}\label{eq:LLie} \NIJ_\Gamma\xrightarrow{\LLie}\NIJ_A \end{equation} 
is a map from the space of multiplicative $(1,1)$-tensors on $\Gamma\toto M$ 
to the space of infinitesimal multiplicative $(1,1)$-tensors on $A$.
The next proposition indicates that this map is indeed a bijection 
when $\Gamma$ is s-connected and s-simply connected.

\begin{prop}\label{prop:inttensors}
Let $\Gamma \toto M$ be a s-connected and s-simply connected
 Lie groupoid with Lie algebroid $(A,\rho,[\cdot,\cdot])$.
The assignment $\LLie$ as in Eq.~\eqref{eq:LLie}
 is a bijection between multiplicative $(1,1)$-tensors
on $\Gamma \toto M $ and infinitesimal multiplicative $(1,1)$-tensors
on $A$.
\end{prop}

\begin{proof}
Let $(\phi_A, \phi_M)$ be an infinitesimal multiplicative $(1,1)$-tensor
on $A$.
Since the  tangent Lie groupoid $T\Gamma \toto TM$ is 
 s-connected and s-simply connected,
there exists a Lie groupoid homomorphism $(\NIJ_\Gamma, \NIJ_M)$:
\[ \xymatrix{
T\Gamma \dar[d] \ar[r]^{{\NIJ_\Gamma} } & T\Gamma \dar[d] \\
TM \ar[r]_{{\NIJ_M}} & TM
} \]
which integrates the Lie algebroid morphism \eqref{eq:Axxx}.
It remains to show that $\NIJ_\Gamma$ is a $(1,1)$-tensor on
$\gm$.

Since $\NIJ_A: TA \to TA $ is a $(1,1)$-tensor, we have  the following
commutative diagram of Lie algebroid morphisms:
\[ \xymatrix{ TA \ar[r]^{\NIJ_A} \ar[d]_{q_A} & TA \ar[d]^{q_A} \\
A \ar[r]_{\idn } & A } \]
which implies the following
commutative diagram of Lie groupoid homomorphisms:
\[ \xymatrix{ T\Gamma \ar[r]^{\NIJ_\Gamma} \ar[d]_{q_\Gamma} & T\Gamma \ar[d]^{q_\Gamma} \\
\Gamma \ar[r]_{\idn} & \Gamma } \]
Hence, $\NIJ_\Gamma: T\Gamma \to T \Gamma $ is a map over the
identity of $\Gamma $. It remains to check that it is fiberwise linear.

We recall some well-known constructions in \cite{MX2}.
%[INTEGRATION OF BIALGEBROIDS, LEMMA 4.3].
Note that there is a natural Lie groupoid structure on
 $ T\Gamma \times_\Gamma T\Gamma \toto TM \times_M TM$,
 and the fiberwise addition map:
\[ +_{\Gamma} (V,W) = V+W \qquad \forall V,W 
\in T_\gamma \Gamma, \gamma \in \Gamma \]
is a Lie groupoid homomorphism from $ T\Gamma \times_\Gamma T\Gamma \toto  TM
  \times_M TM  $ to $T\Gamma \toto TM$.
On the level of Lie algebroids,  there is a natural Lie algebroid
 structure on the vector bundle $ TA \times_A TA \to TM \times_M TM$,
 and the fiberwise addition map:
\[ +_{A} (V,W) = V+W \qquad \forall  V,W \in T_a A, \ a \in A \]
  is a Lie algebroid morphism from $ TA \times_A TA \to  TM  \times_M TM$ to $TA \to TM $. 
  The Lie groupoid  homomorphism $\NIJ_\Gamma $ induces a Lie groupoid 
homomorphism $ \NIJ_\Gamma \times \NIJ_\Gamma $
  from $ T\Gamma \times_\Gamma T\Gamma \toto TM \times_M TM $ to itself, while
  the Lie algebroid morphism $\NIJ_A$ induces a Lie algebroid morphism $ \NIJ_A \times \NIJ_A $
  from $ TA \times_A TA \to TM \times_M TM $ to itself.

Moreover, we have 
\[ \left\{ 
\begin{aligned}
& \LLie(T\Gamma\times_\Gamma T\Gamma)=TA\times_A TA \\
& \LLie(\NIJ_\Gamma\times\NIJ_\Gamma)=\NIJ_A\times\NIJ_A \\ 
& \LLie(+_{\Gamma})=+_A 
\end{aligned} 
\right. \]

Since 
\[ \xymatrix{ TA \times_A TA \ar[r]^{\NIJ_A \times \NIJ_A } \ar[d]_{+_{A}} & TA \times_A TA \ar[d]^{+_A} \\ TA \ar[r]_{\NIJ_A} & TA } \] 
is a commutative diagram of Lie algebroid morphisms, 
\[ \xymatrix{ 
T\Gamma\times_\Gamma T\Gamma \ar[r]^{\NIJ_\Gamma\times\NIJ_\Gamma} \ar[d]_{+_{\Gamma}} & T\Gamma\times_\Gamma T\Gamma \ar[d]^{+_{\Gamma}} \\
T\Gamma \ar[r]_{\NIJ_\Gamma} & T\Gamma } \]
is a commutative diagram of Lie groupoid homomorphisms.
Therefore the map $\NIJ_\Gamma: T\Gamma \to T \Gamma $ respects the
 fiberwise addition of the
vector bundle $q_\gm : T\gm\to \gm$.
 Since $\NIJ_\Gamma$ is a smooth map, it follows
that  $\NIJ_\Gamma$ must be fiberwise linear. Hence, it is 
a multiplicative $(1,1)$-tensor on $\gm$. 

Finally, it is simple to see, from the construction, that
\[ \LLie(\NIJ_\Gamma)=\phi_A .\]
This concludes the proof.
\end{proof}

\subsection{Multiplicative Nijenhuis tensors on Lie groupoids}

We now can state one of the main theorems of this section.

\begin{thm}\label{thm:inttensors}
Let $\Gamma \toto M$ be a s-connected and s-simply connected 
Lie groupoid with Lie algebroid $(A,\rho,[\cdot,\cdot])$.
The assignment $\LLie$ as in Eq.~\eqref{eq:LLie} 
is a bijection between multiplicative Nijenhuis tensors
on $\Gamma\toto M$ and infinitesimal multiplicative Nijenhuis tensors on $A$.
\end{thm}

\begin{proof}
Theorem~\ref{th:Lie_deriv} implies that 
\[ \LLie(\mathcal{N}_{\NIJ_\Gamma})=\mathcal{N}_{\LLie(\NIJ_\Gamma)}=\mathcal{N}_{\NIJ_A} .\]
Since $\NIJ_A:TA\to TA$ is a Nijenhuis tensor, the Nijenhuis torsion 
$\mathcal{N}_{\NIJ_A}$ vanishes. The latter is equivalent to
the commutativity of the diagram
\[ \xymatrix{ \wedge^2TA \ar[r]^{\mathcal{N}_{\NIJ_A}} \ar[d]_{q_A} & TA \\
A \ar[ru]_{\imath_A} & } \]
where $\imath_A : A \hookrightarrow  TA $ is the embedding of
the zero section.
Since all the maps in the above diagram are Lie algebroid morphisms,
it follows from the relations  $q_A=\LLie (q_\gm)$ and $\imath_A
=\LLie (\imath_\gm )$ that
one has a commutative diagram of Lie groupoid homomorphisms
\[ \xymatrix{ \wedge^2T\gm \ar[r]^{\mathcal{N}_{\NIJ_\gm}} \ar[d]_{q_\gm} & T\gm \\
\gm \ar[ru]_{\imath_\gm} & } \]
where $\imath_{\Gamma} : \Gamma  \to T\Gamma $ is the  embedding to
the zero section. 
Therefore the Nijenhuis torsion ${\mathcal N}_{\NIJ_\Gamma} $ vanishes.
This completes the proof of the theorem.
\end{proof}

\subsection{Holomorphic Lie groupoids}
\label{section:7.3}

\begin{defn}: A holomorphic Lie groupoid is a (smooth) Lie groupoid $\Gamma\toto X$,
where both $\Gamma$ and $X$ are
 complex manifolds and  all the
structure maps $m$, $\unit$, $\inverse$, $s$ and $t$ are holomorphic maps.
\end{defn}

This definition requires a justification. The manifold
\[ \Gamma_2 = \{(\gamma_1,\gamma_2 ) | t(\gamma_1) = s(\gamma_2), 
\gamma_1, \gamma_2 \in \gm \} \]
admits a natural complex manifold structure since both the source and
target maps are holomorphic surjective submersions.
It thus makes sense to require the multiplication map to be holomorphic.

The following result essentially follows from the definition and 
the Newlander-Nirenberg theorem.

\begin{prop}\label{prop:7.5}
The following assertions are equivalent:
\begin{enumerate}
\item $\Gamma\toto X$ is a holomorphic Lie groupoid;
\item $\Gamma\toto X$ is a real Lie groupoid, where
 $\Gamma$ and $X$ are almost complex manifolds with associated
 almost complex structures $J_{\Gamma}$ and $J_X$
 respectively such that $(J_{\Gamma}, J_X)$ is a
multiplicative Nijenhuis tensor on $\gm\toto X$.
\end{enumerate}
\end{prop}

A holomorphic Lie algebroid is said to be \emph{integrable} 
if it is the Lie algebroid associated to some holomorphic Lie groupoid.

The main theorem of this section is the following

\begin{thm}
\label{thm:main}
Assume that $A\to X$ is a holomorphic Lie algebroid.
Let $\gm$ be a $s$-connected and $s$-simply connected Lie groupoid
integrating the underlying real Lie algebroid $A_R$.
Then $\gm$ admits a natural complex
structure which makes it into a holomorphic Lie groupoid with
its holomorphic Lie algebroid being $A\to X$.

In other words, a holomorphic Lie algebroid $A$ is integrable
if, and only if, its underlying real Lie algebroid $A_R$
is integrable. Similarly, a holomorphic Lie algebroid is integrable
if, and only if, its underlying imaginary Lie algebroid $A_I$
is integrable.
\end{thm}

\begin{proof} 
According to Proposition~\ref{prop:tgtalg},
the almost complex structure $(J_A,J_X)$ is an infinitesimal
multiplicative Nijenhuis tensor on $A$. By Theorem~\ref{thm:inttensors},
it integrates to a multiplicative Nijenhuis tensor
$(J_{\Gamma},J_X)$ on $\gm\toto X$. 
Since $-\idn_{\Gamma}:T\Gamma\to T\Gamma$
is also a multiplicative $(1,1)$-tensor and $\LLie(-\idn_{\Gamma})=-\idn_A$,
it follows that $J_\Gamma^2=-\idn_{\Gamma}$ and therefore
$J_\Gamma$ is an almost complex structure on $\gm$.
The conclusion thus follows from Proposition~\ref{prop:7.5}.
\end{proof}

As a consequence, we can determine whether a holomorphic
Lie algebroid is integrable by applying 
 the integrability criteria of Crainic-Fernandes:
Theorem 4.1 in \cite{CF} to  its underlying real
Lie algebroid.

\begin{rmk}
Note that, a complex structure on the Lie algebra of a real Lie group G 
extends uniquely to an integrable complex Lie group structure on G. 
This is false for groupoids. The condition of $s$-connectedness and $s$-simply connectedness 
is indeed necessary.
For instance, take a bundle of groups $X\times\CC\toto X$ and consider it as a groupoid,
where $X$ is a complex manifold and $\CC$ is equipped with the additive group structure.
Let $L$ be a smooth non-holomorphic function from $X$ to the space of lattices in $\CC$
(for instance $X=\CC$ and $f(x+iy)=1+isin(x)$) and consider the quotient groupoid 
$(X\times\CC)/L\toto X$. Clearly this is a real
Lie groupoid integrating the underlying real 
Lie algebroid $X\times\CC\to X$ with zero anchor and zero bracket. However the 
holomorphic structure on the Lie algebroid does not extend to a holomorphic groupoid structure 
on this quotient groupoid $(X\times\CC)/L\toto X$.
\end{rmk}

A holomorphic Lie algebroid may not be always integrable, as shown in the following 

\begin{example}
Let $X$ be a complex manifold and $\omega$ a  holomorphic
closed 2-form on $X$. Then $A: TX\oplus (X\times {\mathbb C})\to X$
is naturally equipped with a holomorphic Lie algebroid structure,
where the anchor is the projection onto the
first component, and the Lie bracket is
$$ [(X,f),(Y,g)] = ([X,Y], X(g)-Y(f)+\omega(X,Y)),$$
$\forall X, Y\in \shs{A} (U), \ f, g\in \hf{X}(U)$.
It is simple to see that its  underlying real Lie algebroid $A_R$
is isomorphic to $TX\oplus (X\times {\mathbb R}^2)\to X$,
where the Lie bracket is given by
$$ [(X,f_1, f_2),(Y,g_1, g_2)] = ([X,Y], X(g_1)-Y(f_1)+\omega_1(X,Y), X(g_2)-Y(f_2)+\omega_2(X,Y)),$$
$\forall X, Y\in  {\mathfrak X} (X)$ and $f_1, f_2, g_1, g_2\in C^{\infty}({\mathbb R}^2, {\mathbb R})$, where $\omega_1$ and $\omega_2$ are
the real and imaginary parts of $\omega$, respectively, i.e. 
$\omega=\omega_1 +i\omega_2 $.

According to Example 3.7 in \cite{CF},
 the   Lie algebroid  $TX\oplus (X\times {\mathbb R}^2)\to X$
 is integrable if and only if 
  the group of periods of $(\omega_1, \omega_2)$:
$$\{  (\int_\gamma \omega_1, \int_\gamma \omega_2)| \ \   [\gamma] \in
 \pi_2 (M,x) \} $$
is a discrete subgroup of ${\mathbb R}^2$
(the argument in \cite{CF} is only presented for ${\mathbb R}$-valued  closed
$2$-forms, but  it clearly extends to ${\mathbb R}^2$-valued  closed
 $2$-forms).

%In order to find a counter example, one has therefore to find a
% holomorphic $2$-form $\omega $ for which the group
%of periods of $\omega $, seen as a ${\mathbb C} $-valued smooth
%$2$-form, is not discrete. This can be made as follows.

Consider the complex manifold
 $$N=\left\{  (z_1, z_2 ,z_3) \in {\mathbb C}^3  \left| z_1^2 + z_2^2 + z_3^2 =1
 \right. \right\}$$
and the holomorphic 2-form on ${\mathbb C}^3$:
$$ \eta = \frac{1}{4\pi^2} \big(z_1 \diff z_2 \wedge \diff z_3 + z_2 \diff z_3
\wedge \diff z_1  +  z_3 \diff z_1 \wedge \diff z_2  \big) .$$
It is clear that the pull back of $\eta$
defines   a holomorphic closed 2-form on $N$.
By $\eta_1$ and $\eta_2$,  we denote the real part and 
the imaginary part of $\eta$, respectively.
Consider  the submanifold ${\mathbb R}^3$ of ${\mathbb C}^3$
defined by ${\mathbb R}^3\cong \{(z_1, z_2, z_3)|y_1=y_2=y_3=0\}$.
The intersection $ S$ of $N$ with ${\mathbb R}^3$ is a unit sphere:
 $$S = N \cap {\mathbb R}^3=\left\{  (x_1,x_2 ,x_3) \in {\mathbb R}^3  \left|
x_1^2 + x_2^2 + x_3^2 =1 \right. \right\}  $$
The pull back of $\eta_1$ to $S$ is a real valued
$2$-form, which takes the same expression:
$$ \eta_1|_{S} = 
\frac{1}{4\pi^2} \big(x_1 \diff x_2 \wedge \diff x_3 + x_2 \diff x_3
\wedge \diff x_1  +  x_3 \diff x_1 \wedge \diff x_2  \big) .$$
It is clear that $\eta_1$ is a volume form on the unit sphere:
 $\int_S \eta_1 = 1$. On  the other hand, the pull back of $\eta_2$
to $S$ vanishes.

Let $X = N \times N $, and $\omega = \sqrt{2} \,  p_1^* \eta +
p_2^* \eta $, where $ p_1,p_2$ stand for the projections on the first
and second components.
Then $\omega$ is a holomorphic closed 2-form on $X$.
It is simple to see that the group of periods  of $(\omega_1, \omega_2)$
 contains $(1, 0)$ and $(\sqrt{2}, 0)$, and therefore is not discrete.
%: the first value is
%obtained by integrating $\omega $ along $S \times \{a\} $, with $a \in
%N$ a fixed point, while the second one is obtained  by integrating
% $\omega $ along $\{a\} \times S $.
Hence, the corresponding holomorphic Lie algebroid $A$ is not integrable.  
\end{example}

We end this section with the following

\begin{prop}\label{prop:hol_mor}
Let $\Gamma \toto X $ and $ \Gamma' \toto X'$
be holomorphic Lie groupoids with holomorphic Lie algebroids
$A\to X$ and $A'\to X'$, respectively. 
Assume that $\phi: \Gamma \to \Gamma'$ is a  homomorphism
of the underlying real Lie groupoids
and its infinitesimal $\Lie (\phi): A\to A'$ is a morphism of
 holomorphic Lie algebroids.
Then $\phi$ is a homomorphism of holomorphic Lie groupoids.
\end{prop}

\begin{proof}
The map $ \Lie (\phi): A\to A'$ is holomorphic implies that the following 
\[ \xymatrix{ T A \ar[r]^{T\Lie (\phi)}\ar[d]_{J_{A}} & T A'  \ar[d]^{J_{A'}} \\
 T A \ar[r]_{T\Lie (\phi)} &  T A' } \]
is a commutative diagram of Lie algebroid morphisms.

Hence it implies the commutativity of the
Lie groupoid homomorphisms of the diagram below:
\[ \xymatrix{ T \Gamma \ar[r]^{T\phi}\ar[d]_{J_{\Gamma}} & T \Gamma'  \ar[d]^{J_{\Gamma'}} \\
 T \Gamma \ar[r]_{T\phi} &  T \Gamma' .} \]
The latter exactly means that $\phi: \Gamma \to \Gamma' $
is a holomorphic map.
\end{proof}

\subsection{Holomorphic symplectic groupoids}

As an application, in this section,
we study holomorphic symplectic groupoids.

Holomorphic symplectic groupoids can be introduced in
a similar fashion as in the smooth case.

\begin{defn}
A holomorphic symplectic groupoid is a holomorphic
Lie groupoid $\gm\toto X$ together with a holomorphic symplectic
2-form $\omega\in\Omega^{2,0}(\gm)$ such that
the graph of multiplication $\Lambda\subset\gm\times\gm\times\bar{\gm}$ 
is a Lagrangian submanifold, where $\bar{\gm}$
stands for the $\gm$ equipped with the opposite symplectic structure.
\end{defn}

As in the smooth case, this last condition is equivalent to
$\partial\omega =0$ where $\partial:\Omega^2(\Gamma)\to
\Omega^2(\Gamma_2)$ is the alternate sum of the pull back maps 
of the three face maps $\gm_2\to\gm$.
In this case, $\omega$ is said to be \emph{multiplicative}.

As in the smooth case, if $\gm\toto X$ is a holomorphic
symplectic groupoid, then $X$ is naturally a holomorphic
Poisson manifold. More precisely, there exists a unique
holomorphic Poisson structure on $X$ such that the source map
$s:\Gamma\to X$ is a holomorphic Poisson map, while
the target map is then an anti-Poisson map. 

Given a holomorphic symplectic groupoid $\gm\toto X$, 
its holomorphic Lie algebroid is isomorphic to the cotangent
Lie algebroid $(T^*X)_\pi\to X$, where $\pi$ is the induced 
holomorphic Poisson structure on $X$.

Conversely, a holomorphic Poisson manifold $(X,\pi)$ is said 
to be \emph{integrable} if it is the induced
holomorphic Poisson structure on the unit space of 
a holomorphic symplectic groupoid $\gm\toto X$.
We say that $\gm\toto X$ {\em integrates} 
the holomorphic Poisson structure $(X,\pi)$.

The main theorem is the following:

\begin{thm}
\label{thm:5.21}
A holomorphic Poisson manifold is integrable if, and only if, 
either its real or its imaginary part is integrable as a real Poisson manifold. 

More precisely, if $(\gm\toto X, \omega_R +i \omega_I)$ is a holomorphic
symplectic groupoid integrating the holomorphic Poisson
structure $(X, \pi_R +i \pi_I)$, then
$(\gm\toto X,  {4}\omega_R)$ and $(\gm\toto X, - {4} \omega_I)$
 are symplectic groupoids integrating  the real Poisson manifolds
 $(X, \pi_R)$ and $(X, \pi_I)$, respectively.   

Conversely given a holomorphic Poisson manifold $(X,\pi)$,
where $\pi=\pire+i \piim\in \gm (\wedge^2 T^{1, 0}X)$,
if $(\gm\toto X, \omega_R)$ is a $s$-connected and $s$-simply connected symplectic groupoid integrating $(X,\pire)$, then 
\begin{enumerate}
\item $\gm\toto X$ admits  a holomorphic Lie groupoid structure. By
  $J_\gm:T\gm\to T\gm$  we
denote its almost complex structure;
\item $(\gm\toto M,\omega_I)$, where $\omega_I(\cdot,\cdot):=\omega_R(J_\Gamma\cdot,\cdot)$,
 is a symplectic groupoid integrating $(X,  \piim )$;
 \item $(\Gamma\toto X,\omega)$, where $\omega:= \frac{1}{4}
(\omega_R-i\omega_I)$ ,
 is a holomorphic symplectic groupoid integrating $(X, \pi)$. 
\end{enumerate}
\end{thm}

This theorem can be derived from Theorem~5.2 in \cite{PqN}
or Theorem~3.4 in \cite{crainic} via
 the equivalence relation between holomorphic Poisson
manifolds and Poisson Nijenhuis structures established
in Theorem ~\ref{PNGC}. We, however, give a direct proof below, as an application of Theorem~\ref{thm:5.21}. We start with a
simple fact of complex geometry.

\begin{lem}
\label{lem:5.22}
Let $X$ be a complex manifold with almost complex structure
$J$ and $\alpha_R \in \Omega^2 (X)$ a two-form
on the underlying real manifold. By $\alpha_R^b$ we denote
its induced bundle map $ \alpha_R^b: TX\to T^*X$.
\begin{enumerate}
\item If $J^*\rond \alpha_R^b=\alpha_R^b \rond J$,
then $\alpha(\cdot , \cdot ):=\alpha (\cdot , \cdot )-i\alpha (J \cdot , \cdot )$ is a $(2, 0)$-form on $X$.
\item If, moreover, $\alpha_R^b: TX\to T^*X$ is a holomorphic map, then
$\alpha$ is a holomorphic 2-form.
\item Furthermore, $d\alpha_R=0$ iff $\partial\alpha=0$.
\end{enumerate}
\end{lem}

\begin{proof}[Proof of Theorem~\ref{thm:5.21}] 
Assume that  $(\gm\toto X, \omega_R +i \omega_I)$ is a holomorphic
symplectic groupoid integrating the holomorphic Poisson
structure $(X, \pi_R +i \pi_I)$. Clearly, both
 $ (\gm \toto X, \omega_R)$ and $ (\gm \toto X, \omega_I)$
are symplectic  groupoids.
By $\Pi_R$ and $\Pi_I$ we denote the Poisson tensors corresponding to
the symplectic two-forms $\omega_R$ and $\omega_I$, respectively.
%By Proposition~\ref{prop:2.9},
Then the holomorphic Poisson
bivector field corresponding to $\omega$ is 
$\frac{1}{4}(\Pi_R-i \Pi_I)$.
A Poisson map between two holomorphic Poisson manifolds
is also a Poisson map  between their real parts,  and  as well as
between their imaginary parts.  As a consequence,  we have
$ \frac{1}{4}s_*\Pi_R=\pi_R$ and $-\frac{1}{4}s_*\Pi_I=\pi_I$.
It follows  that
$(\gm \toto X, 4\omega_R)$ and $ (\gm \toto X, -4 \omega_I)$
are symplectic groupoids integrating, respectively, the real part $\pi_R $
and the imaginary part $\pi_I $ of $\pi$.

Conversely, let $(X,\pi)$, where $\pi=\pire+i\piim\in\gm(\wedge^2 T^{1,0}X)$, 
be a holomorphic Poisson manifold and $(\gm\toto X,\omega_R)$ 
a $s$-connected and $s$-simply connected symplectic groupoid integrating $(X,\pire)$.
Hence the Lie algebroid of $\gm\toto X$ is
isomorphic to $(T^*X)_{\pire}$. By Proposition~\ref{prop:3.7}, 
$(T^*X)_{\pire}$ is the underlying real Lie algebroid
of the holomorphic Lie algebroid $(T^*X)_{\frac{1}{4}\pi}$.
According to Theorem~\ref{thm:main}, $\gm\toto X$
admits a multiplicative almost complex structure
$J_\gm$ which makes $\gm$ into a holomorphic Lie groupoid.

For the remaining of the proof, we need recall some
well known facts concerning Poisson groupoids \cite{MX, MX2}.
For a Poisson groupoid $(\gm\toto X, \Pi)$, the map
\begin{equation}\label{eq:53}\xymatrix{
 T^*\Gamma \dar[d] \ar[r]^{ \Pi^{\#} } & T\Gamma \dar[d] \\
 A^* \ar[r]_{\rho_* } & TX }\end{equation}
induced by the Poisson tensor $\Pi$ is a Lie groupoid homomorphism,
where $(A,A^*)$ denotes its corresponding Lie
bialgebroid with anchors $\rho$ and $\rho_*$, respectively.
Here $T\gm\toto TX$ and $T^*\gm\toto A^*$ are, respectively, 
the tangent and cotangent Lie groupoids of $\gm\toto X$ \cite{CDW}.
The infinitesimal $\Lie(\Pi)$ of Eq.~\eqref{eq:53} 
is the Lie algebroid homomorphism below (Theorem~8.3 in \cite{MX}):
\begin{equation*}
\xymatrix{ T^*A \ar[r]^{\pi_{A}^{\#} }\ar[d]_{} & TA \ar[d]^{} \\
A^* \ar[r]_{a_*} & TX .}
\end{equation*} 
where the left hand side $T^{*}A\lon A^{*}$ is the
cotangent Lie algebroid $T^{*}(A^{*})\lon A^{*}$ 
induced by the Lie Poisson structure on $A^{*}$ by identifying
$T^{*}A$ with $T^{*}(A^{*})$ via the canonical isomorphism 
$R\colon T^{*}A^{*}\lon T^{*}A$ as given in Theorem~5.5 of \cite{MX}, 
the right hand side $TA\lon TX$ is the tangent Lie algebroid of $A$, 
and $\pi_{A}$ is the Lie Poisson tensor on $A$.

Now consider the symplectic
groupoid $(\gm\toto X, \omega_R)$. Its corresponding
Lie bialgebroid is $((T^*X)_{\pi_R}, TX)$,
and the Lie Poisson structure on $T^*X$
dual to the Lie algebroid $TX$ is the canonical
symplectic structure on  the cotangent bundle. Let
$\alpha \in \Omega^2 (T^*X)$ be its symplectic form.
Applying the above general result on Poisson
groupoids and reversing the maps, we have
a Lie groupoid homomorphism
\begin{equation} \label{eq:55} \xymatrix{ 
T\Gamma \dar[d] \ar[r]^{ \omega_R^b } & T^*\Gamma \dar[d] \\ 
TX \ar[r]_{\idn } & TX } \end{equation}

Its infinitesimal $\Lie(\omega_R^b)$ is the Lie algebroid homomorphism:
\begin{equation} \label{eq:56} \xymatrix{ 
T(T^*X) \ar[r]^{\alpha^b}\ar[d]_{} & T^* (T^*X) \ar[d]^{} \\
TX \ar[r]_{\idn} & TX .} \end{equation}

So far we have ignored the complex structures, and all
the maps involved in Eqs.~\eqref{eq:55} and \eqref{eq:56} are considered as maps of real manifolds. 
On the other hand, $T\gm \toto TX$ and $T^*\gm \toto TX$ are 
both holomorphic Lie groupoids since $\gm\toto X$ is a
holomorphic Lie groupoid. And their corresponding
holomorphic Lie algebroids $T^*(T^*X)\to TX$ and $T (T^*X)\to TX$ are, respectively, the holomorphic
counterparts of the Lie algebroids involved in Eq.~\eqref{eq:56}. 
It is easy to see that 
$\alpha^b$ is a holomorphic map, and therefore is
 indeed a holomorphic Lie algebroid homomorphism.
By Proposition~\ref{prop:hol_mor}, it follows that $\omega_R^b$ is
indeed a homomorphism
of holomorphic Lie groupoids. Moreover, the diagram 
\begin{equation*}
\xymatrix{ T(T^*X) \ar[r]^{\alpha^b}\ar[d]_{J_{T^*X}} & T^* (T^*X) \ar[d]^{J_{T^*X}^*} \\
T(T^*X) \ar[r]_{\alpha^b} & T^*(T^*X) }
\end{equation*}
is a commutative diagram of real Lie algebroid homomorphisms, since
$\alpha$ is essentially the real part of the canonical holomorphic
symplectic 2-form on $T^*X$. Hence we
have a commutative diagram of Lie groupoid homomorphisms:
\begin{equation*}
\xymatrix{ T\gm \ar[r]^{\omega_R^b}\ar[d]_{J_{\gm}} & T^*\gm \ar[d]^{J^*_{\gm}} \\ T\gm \ar[r]_{\omega_R^b} & T^*\gm .}
\end{equation*}

Therefore, by Lemma~\ref{lem:5.22}, we conclude that
$\omega=\omega_R-i\omega_I\in\Omega^{2,0}(\gm)$,
where $\omega_I(\cdot,\cdot)=\omega(J\cdot,\cdot)$, 
is a holomorphic closed 2-form. It is clear that
$\omega$ is non-degenerate and multiplicative. 
Hence $(\gm\toto X,\omega)$ is
a holomorphic symplectic groupoid. Finally we have
\[ s_*(\omega_R^b J_\gm)^{-1}s^*=-s_* J_\gm(\omega_R^b)^{-1}s^*
=-J_X s_*(\omega_R^b)^{-1}s^*=-J_X\pi_R^{\#}=\pi_I^{\#}. \]
The conclusion thus follows. 
%then $(X,\pi_I, J) $ is a Poisson Nijenhuis structure.
% According to \cite{PqN},
%$(\pi_I, J )$ integrates to a symplectic Nijenhuis 
%Lie groupoid $(\Gamma \toto M, \omega_I, J_\Gamma)$.
%One easily sees that $J_\Gamma^2=-\idn$. Hence $\Gamma \toto M$
%is a holomorphic  symplectic groupoid.
\end{proof}

Indeed the exact same proof as above leads to the following

\begin{thm}
Let $(X,\pi)$, where $\pi=\pire+i\piim\in\gm(\wedge^2 T^{1,0}X)$, 
be a holomorphic Poisson manifold.
Assume that $\gm\toto X$ is a holomorphic Lie groupoid with
almost complex structure $J_\gm$, whose corresponding holomorphic Lie algebroid is $(T^*X)_{\tfrac{1}{4}\pi}$.
Moreover, assume that there exists a symplectic
real 2-form $\omega_R$ on the underlying real Lie groupoid
such that $(\gm\toto X,\omega_R)$ is a symplectic groupoid
integrating the real Poisson structure $\pire$. Then
$\omega=\omega_R-i\omega_I\in\Omega^{2,0}(\gm)$,
where $\omega_I(\cdot,\cdot)=\omega(J\cdot,\cdot)$, 
is a multiplicative holomorphic symplectic 2-form
on $\gm$ so that $(\gm\toto X,\frac{1}{4}\omega)$ is
a holomorphic symplectic groupoid integrating
the holomorphic Poisson structure $\pi$.
\end{thm}

\subsection{Holomorphic extension of analytic Poisson structures}

To any analytic Poisson structure $\pi_{real}$ on ${\mathbb R}^n$,
there associates a canonical  holomorphic  Poisson structure $\pi$
on ${\mathbb C}^n$, called its  {\em holomorphic extension}, which
 can be defined as follows. Assume that the Poisson 
brackets of $\pi_{real}$, on the canonical coordinates,
 are given by
$$\{x_i, x_j\}=\phi_{ij}(x_1, \dots, x_n ),$$
where  $\phi_{ij}$ are, for all $i,j=1 , \dots, n$,   real
 analytic functions on ${\mathbb R}^n$.
Then the Poisson 
brackets of its holomorphic extension on ${\mathbb C}^n$ are
given by $$\{z_i, z_j\}=\phi_{ij}(z_1, \cdots, z_n ).$$
 
Then we have the following:

\begin{thm}\label{thm:integrability_extension}
If the holomorphic extension $\pi$ of a real analytic Poisson
structure $\pi_{real} $ is integrable, then $\pi_{real} $ 
must be  integrable.
\end{thm}

Let $\pi_R$ and $\pi_I $ be the real part and the imaginary part of $\pi$. We advise the reader to keep in mind that
$\pi_R $ is a smooth bivector field on ${\mathbb C}^n$,
and not to confuse it with $\pi_{real} $, which is a smooth bivector
field on ${\mathbb R}^n$.

Recall that a Poisson involution on a Poisson manifold $P$
is a Poisson diffeomorphism $\Phi : P\lon P$ such that $\Phi^2 =id$.

\begin{lem}
\label{lem:Poi_invol}
The complex conjugation map $\sigma(z_1, \cdots, z_n) =
 (\overline{z_1}, \cdots, \overline{z_n}) $ is a Poisson involution of 
$({\mathbb C}^n , \pi_R)$.
\end{lem}

\begin{proof}
We write $\phi_{ij}= f_{ij}+ i g_{ij}$, where $f_{ij},\ g_{ij}
\in C^\infty ({\mathbb C}^n, {\mathbb R})$, and $z_k=x_k+i y_k$, $k=1, \cdots, n$. 
One computes immediately  that
\begin{eqnarray}
 \begin{array}{rcl} \pi_R &=&
 \frac{1}{4}\sum_{i,j=1}^n  f_{ij}(z_1,\dots, z_n)
(\frac{\partial}{\partial x_i} \wedge \frac{\partial}{\partial x_j} - \frac{\partial}{\partial y_i} \wedge
 \frac{\partial}{\partial y_j}) \nonumber \\
 & & +  \frac{1}{4}  \sum_{i,j=1}^n  g_{ij}(z_1,\dots, z_n)
 (\frac{\partial}{\partial x_i} \wedge \frac{\partial}{\partial y_j} + \frac{\partial}{\partial y_i} \wedge
 \frac{\partial}{\partial x_j}). \label{eq:piR}
\end{array} 
\end{eqnarray}

By constructions, we have $\sigma^* \phi_{ij}=\overline{\phi_{ij}}$. Thus it
follows that $\forall i, j=1, \cdots, n$,
\begin{equation}
\sigma^* f_{ij} = f_{ij}, \ \ \  \sigma^* g_{ij} = -g_{ij}.
\end{equation}
On the other hand, it is obvious that
\begin{equation}
 \sigma_*  (\frac{\partial}{\partial x_j})
=\frac{\partial}{\partial x_j}  , \ \   \sigma_*  (\frac{\partial}{\partial
  y_j}) =- \frac{\partial}{\partial y_j}
\end{equation}
The conclusion thus follows immediately.
\end{proof}
 
\begin{remark}
Note that $\sigma $ is an anti-Poisson map with respect to the
imaginary part  $\pi_I $ of~$\pi $.
\end{remark}

It is well known that the  stable  locus of a Poisson involution
carries a natural Poisson structure \cite{FV, Xu:03}.
More precisely, 
let $Q$ be the stable locus of a Poisson involution
$\Phi : P\lon P$. Assume that the Poisson tensor $\pi$
on $P$ is $\pi =\sum_{i} X_{i}\wedge Y_{i}$, where $X_{i}$ and $Y_{i}$
are vector fields on $P$. Then the induced 
Poisson tensor $\pi_{Q}$ on $Q$ is given by
\begin{equation}
\label{eq:piQ}
\pi_{Q} =\sum_{i} X_{i}^{+}\wedge Y_{i}^{+} |_{Q},
\end{equation}
where $X_{i}^{+}$ and $ Y_{i}^{+}$ are vector fields
on $P$ defined by $X_i^{+} =\half (X_i+\Phi_* X_i)$ and
$Y_i^{+} =\half (Y_i+\Phi_* Y_i)$, respectively.
The Poisson manifold $(Q, \pi_Q)$ is called
a Dirac submanifold of $(P, \pi)$ (which is also
called a Lie-Dirac submanifold in \cite{CF2}).

The following result was proved in \cite{Xu:03} (see also  \cite{Fer_fp, CF2}).

\begin{prop}
\label{cor:inv-poid}
If $Q$ is the  stable locus of a
Poisson involution on an  integrable Poisson manifold $P$,
    then $Q$ is always an  integrable Poisson manifold itself.
\end{prop}

\begin{lem}\label{lem:lie_dirac}
Let $\pi$ be the holomorphic extension on ${\mathbb C}^n$
of an analytic Poisson structure $\pi_{real}$ on ${\mathbb R}^n$.
The induced Poisson structure  on the stable locus of
the Poisson involution $\sigma$ on 
$({\mathbb C}^n , \pi_R)$ is isomorphic to 
$({\mathbb R}^n,   \frac{1}{4}\pi_{real})$.
\end{lem}
\begin{proof}
Let $Q$ be the stable locus of $\sigma$.
Then $Q=\{(z_1, \cdots, z_n)|y_1=\cdots =y_n=0\}\cong {\mathbb R}^n$.
By Eqs. \eqref{eq:piR} and \eqref{eq:piQ}, we have
$$ (\pi_R)_Q=\frac{1}{4}\sum_{i,j=1}^n  f_{ij}(z_1,\dots, z_n)|_Q
\frac{\partial}{\partial x_i} \wedge \frac{\partial}{\partial x_j}.
 $$
Since $f_{ij}(z_1,\cdots, z_n)|_Q =\phi_{ij} (x_1, \cdots, x_n)$,
the right hand side is equal to $ \frac{1}{4} \pi_{real}$.
\end{proof}

\begin{proof}[Proof of Theorem~\ref{thm:integrability_extension}]
It  follows immediately from Theorem \ref{thm:5.21}, Proposition \ref{cor:inv-poid} and Lemma \ref{lem:lie_dirac}.
\end{proof}

\begin{example}
Consider the holomorphic Poisson structure 
\begin{equation}
\label{eq:pih}
 \pi = e^{\frac{z_1^2 + z^2_2 + z_3^2}{2} } z_1  \frac{\partial}{\partial z_2} \wedge \frac{\partial}{\partial z_3}  + \mbox{ c.p. } 
\end{equation} 
on ${\mathbb C}^3$.
It is clear that this is a holomorphic extension of the
analytic Poisson structure:
\begin{equation}
\label{eq:pireal}
  \pi_{real} = e^{\frac{x_1^2 + x^2_2 + x_3^2}{2} } x_1  \frac{\partial}{\partial x_2} \wedge
 \frac{\partial}{\partial x_3}  + \mbox{ c.p. } 
\end{equation}
on ${\mathbb R}^3$.

According to  Example 3.3 in \cite{CF-lec}, the Poisson
 structure \eqref{eq:pireal} is not integrable. By  Theorem
\ref{thm:integrability_extension}, the holomorphic Poisson
structure \eqref{eq:pih} on ${\mathbb C}^3$ is not integrable either.
\end{example}


\begin{thebibliography}{10}
\bibitem{Boyom}
M. N. Boyom,
\newblock KV-cohomology of Koszul-Vinberg algebroids and Poisson manifolds,
\newblock {\em Internat. J. Math.}{\bf 16} (2005), no. 9, 1033--1061.

                                                       
%\bibitem{Block}
% J. Block,
%\newblock Duality and equivalence of module categories 
%in noncommutative geometry.
%\newblock  arXiv:math/0509284.

%\bibitem{Bondal}
%A. Bondal,
%\newblock Non-commutative deformations and Poisson brackets on
% projective spaces.
%\newblock  Preprint no. 67, Max-Planck-Institut, Bonn 1993.

%\bibitem{BZ}
%J.-L. Brylinski and G.  Zuckerman, 
%\newblock  The outer derivation of a complex Poisson manifold,
%\newblock  {\em J. Reine Angew. Math.} {\bf 506} (1999), 181-189.


\bibitem{Canez}
S. Canez,
\newblock Private communication.

\bibitem{CDW}
A.~Coste, P.~Dazord, and A.~Weinstein, \emph{Groupo\"\i des symplectiques},
  Publications du D\'epartement de Math\'ematiques. Nouvelle S\'erie. A, Vol.\
  2, Publ. D\'ep. Math. Nouvelle S\'er. A, {\bf 87}, Univ. Claude-Bernard, Lyon,
  1987, pp.~i--ii, 1--62.

%\bibitem{Cav}
%Gil R. Cavalcanti.
%\newblock  New aspects of the ddc-lemma
%\newblock  arXiv:math/0501406.

\bibitem{CattaneoF}
  A. Cattaneo and  G. Felder, 
\newblock     Poisson sigma models and symplectic groupoids,
\newblock {\em Prog. Math.} {\bf 198} (2001) 61--93,

\bibitem{crainic}
M. Crainic,
\newblock Generalized complex structures and Lie brackets
\newblock arXiv:math/0412097 

\bibitem{CF} 
M. Crainic and R.-L. Fernandes,
\newblock  Integrability of Lie brackets. 
\newblock  {\em Ann. of Math.} (2)  {\bf 157 } (2003),  575--620.

\bibitem{CF2}
 M. Crainic and R.-L. Fernandes, 
\newblock Integrability of Poisson brackets,
\newblock  {\em J. Differential Geometry}. {\bf 66} (2004),  71--137.

\bibitem{CF-lec}
 M. Crainic and   R.L. Fernandes,
\newblock    Lectures on integrability of Lie brackets,
\newblock  arXiv:0611259v1


%\bibitem{DF}
% P. A. Damianou, R. L.  Fernandes,
% \newblock  Integrable hierarchies and the modular class
%\newblock arXiv:math/0607784


%\bibitem{Drinfeld83a}
%V.G. Drinfel'd, 
%\newblock  Hamiltonian structures on Lie groups, Lie bialgebras,
%and the geometric meaning of the classical Yang-Baxter equations,
%\newblock  {\em Soviet Math. Dokl.} {\bf 27}, (1983), 68-71.

%\bibitem{Drinfeld83b}
%V.G. Drinfel'd, 
%\newblock   On constant quasiclassical solutions of the Yang-Baxter
%quantum equation,
%\newblock {\em Soviet Math. Dokl.} {\bf 28}, (1983), 667-671.

\bibitem{ELW}
S. Evens, J.-H. Lu,  and A. Weinstein, 
\newblock  Transverse measures, the modular class, and a cohomology pairing for Lie algebroids,
\newblock   {\em  Quart. J. Math. Oxford (2)} {\bf 50} (1999), 417-436.
\newblock arXiv:dg-ga/9610008v1

%\bibitem{Cal}
%M.~Gualtieri,
%\newblock  {{Generalized complex geometry}},
%\newblock   \mbox{arXiv:math.DG/0401221}.

\bibitem{Fer_fp} R.L. Fernandes,
\newblock   A note on proper Poisson actions,
\newblock  arXiv:0503147.
   
\bibitem{FV}
Fernandes, R., and Vanhaecke, P.,
\newblock  Hyperelliptic Prym varieties and integrable systems,
\newblock {\em Commun. Math. Phys.}{\bf  221} (2001) 169-196. 
 

%\bibitem{Hitchin:03}
%N.~Hitchin. 
%\newblock{Generalized {C}alabi-{Y}au manifolds}.
%\newblock {\em  Q. J. Math.}
%  \textbf{54} (2003), no.~3, 281--308.

%\bibitem{Hitchin:06}
% N. Hitchin,
%\newblock Instantons, Poisson structures and generalized Kaehler geometry. 
%\newblock {\em Commun. Math. Phys.},  265:131-164, (2006).

%\bibitem{Hue1}
%J. Huebschmann,
%\newblock Poisson geometry of certain moduli spaces,
%\newblock The Proceedings of the Winter School "Geometry and Physics" 
%\newblock {\em Rend. Circ. Mat. Palermo (2)} Suppl. No. 39 (1996), 15--35.
 
\bibitem{Hue2}
J. Huebschmann,
\newblock Duality for Lie-Rinehart algebras and the modular class,
\newblock {\em J. Reine Angew. Math.}{\bf 510} (1999), 103--159.

%\bibitem{huybrechts}
%D. Huybrechts,
%\newblock {\em Complex geometry}.
%\newblock Universitext. Springer-Verlag, Berlin, 2005.
%%\newblock An introduction.

%\bibitem{Iwa}
%K. Iwasaki,
%\newblock Fuchsian moduli on a Riemann surface---its
% Poisson structure and Poincar\'e-Lefschetz duality,
%\newblock {\em Pacific J. Math.}{\bf 155} (1992), no. 2, 319--340.

%\bibitem{Kapustin}
%A. Kapustin,
%\newblock Topological strings on noncommutative manifolds,
%\newblock {\em Int. J. Geom. Methods Mod. Phys.} {\bf 1} (2004), no. 1-2, 49--81.

%\bibitem{Kor}
%D. Korotkin and H. Samtleben, 
%\newblock On the quantization of isomonodromic deformations on the torus,
%\newblock {\em Internat. J. Modern Phys. A}{\bf 12}  (1997), no. 11, 2013--2029. 

\bibitem{yks}
Y.  Kosmann-Schwarzbach,
\newblock The {L}ie bialgebroid of a {P}oisson-{N}ijenhuis manifold.
\newblock {\em Lett. Math. Phys.}, 38(4):421--428, 1996.

\bibitem{PN}
Y.  Kosmann-Schwarzbach and F.  Magri,
\newblock Poisson-{N}ijenhuis structures.
\newblock {\em Ann. Inst. H. Poincar\'e Phys. Th\'eor.}, 53(1):35--81, 1990.

%\bibitem{KSM1}
%Y.  Kosmann-Schwarzbach and F.  Magri,
%\newblock On the modular classes of Poisson-Nijenhuis manifolds,
%\newblock arXiv:math/0611202.

%\bibitem{LWX}
%Z.-J. Liu, A. Weinstein, and P. Xu,
%\newblock Manin triples for {L}ie bialgebroids.
%\newblock {\em J. Differential Geom.}, 45(3):547--574, 1997.

\bibitem{LSX}
C.-G. Laurent, M. Stienon, and P. Xu,
Holomorphic Poisson  manifolds and holomorphic Lie algebroids,
\newblock Preprint 2008.
%\bibitem{lu}
%J.-H. Lu,
%\newblock Poisson homogeneous spaces and {L}ie algebroids associated to
%  {P}oisson actions.
%\newblock {\em Duke Math. J.}, 86(2):261--304, 1997.

%\bibitem{mackenzie}
%K.~C.~H. Mackenzie,
%\newblock Ehresmann doubles and drinfel'd doubles for Lie algebroids and Lie
%  bialgebroids.
%\newblock  arXiv:math/0611799.

\bibitem{MX}
K.~C.~H. Mackenzie and P.  Xu,
\newblock Lie bialgebroids and {P}oisson groupoids.
\newblock {\em Duke Math. J.}, 73(2):415--452, 1994.

\bibitem{MX2}
K.~C.~H. Mackenzie and P.  Xu,
\newblock Integration of {L}ie bialgebroids.
\newblock {\em Topology}, 39(3):445--467, 2000.

\bibitem{MR773513}
F.~Magri and C.~Morosi,
 \newblock {On the reduction theory of the {N}ijenhuis
  operators and its applications to {G}elfand-{D}iki\u\i\ equations},
\newblock  Proceedings of the IUTAM-ISIMM symposium on modern
 developments in analytical mechanics, Vol. II (Torino, 1982), {\bf 117}, 1983, pp.~599--626.

\bibitem{MR900387}
F.~Magri and C.~Morosi,
\newblock  \emph{Old and new results on recursion operators: an algebraic
  approach to {KP} equation},
\newblock  Topics in soliton theory and exactly solvable
  nonlinear equations (Oberwolfach, 1986), World Sci. Publishing, Singapore,
  1987, pp.~78--96.

%\bibitem{Mor}
%R. Moraru,  
%\newblock  Integrable systems associated to a Hopf surface.
%\newblock {\em Canad. J. Math.}{\bf 55} (2003), no. 3, 609--635.

%\bibitem{mokri}
%T. Mokri,
%\newblock Matched pairs of {L}ie algebroids.
%\newblock {\em Glasgow Math. J.}, 39(2):167--181, 1997.

%\bibitem{NewNir}
%A.~Newlander and L.~Nirenberg.
%\newblock Complex analytic coordinates in almost complex manifolds.
%\newblock {\em Ann. of Math. (2)}, 65:391--404, 1957.

%\bibitem{Polishchuk}
%A.  Polishchuk,
%\newblock Algebraic geometry of Poisson brackets. Algebraic geometry, 7
%\newblock {\em  J. Math. Sci.} (New York) 84 (1997),  1413--1444. 

%\bibitem{STS}
%M.~A. Semenov-Tian-Shansky, 
%\newblock  Dressing transformations and Poisson Lie group actions,
%\newblock {\em Publ. RIMS, Kyoto University} {\bf 21} (1985), 1237-1260.

\bibitem{PqN}
M. Sti{\'e}non and P. Xu,
\newblock Poisson quasi-{N}ijenhuis manifolds.
\newblock {\em Comm. Math. Phys.} {\bf 27}  (2007),  709--725.

%\bibitem{Sto}
%L. Stolovitch,
%\newblock  Sur les structures de Poisson singuli\`eres,
%\newblock {\em Ergodic Theory Dynam. Systems}{\bf 24} (2004), no. 5, 1833--1863.

%\bibitem{MR1390832}
%I.~Vaisman,
%\newblock  {Complementary {$2$}-forms of {P}oisson structures},
%\newblock  {\em Compositio Math.} \textbf{101} (1996), no.~1, 55--75.

\bibitem{Weinstein}
A.  Weinstein,
\newblock  The integration problem for complex Lie algebroids, 
\newblock {\em From Geometry to Quantum Mechanics, in Honor of Hideki Omori},
Y. Maeda, T. Ochiai, P. Michor, and A. Yoshioka, eds., Progress in
Mathematics, Birkh\"auser, New York (2007), 93-109.

%\bibitem{Xu}
%P.~Xu, {Gerstenhaber algebras and {BV}-algebras in {P}oisson geometry},
%\newblock {\em  Comm. Math. Phys.}
%\newblock  \textbf{200} (1999), no.~3, 545--560.

\bibitem{Xu:03}
P. Xu,
\newblock  Dirac submanifolds and Poisson involutions,
\newblock {\em Ann. Scient. Ec.  Norm.  Sup.} {\bf 36} (2003), 403--430.


\end{thebibliography}
\end{document}